\documentclass[12pt]{article}
\usepackage{amsmath,array,amssymb,latexsym}
\newtheorem{theorem}{Theorem}[section]
\newtheorem{proposition}[theorem]{Proposition}
 
\newtheorem{lemma}[theorem]{Lemma}

\newtheorem{claim}[theorem]{Claim}
\newtheorem{definition}[theorem]{Definition}
\numberwithin{equation}{section}
\begin{document}
\begin{center}
\Large{Lagrangian Mean Curvature Flows and Moment maps}
\\ \vspace{3mm}
\large{Hiroshi Konno\footnote{Supported in part by JSPS KAKENHI Grant Number (C) 23540072 and JP17K05231 \\ Department of Mathematics, School of Science and Technology, Meiji University, \\ 1-1-1 Higashi-Mita, Tama-ku, Kawasaki, Kanagawa 214-8571 Japan
\\
email: hkonno@meiji.ac.jp}}
\end{center}
\begin{abstract}
In this paper, we construct various examples of Lagrangian mean curvature flows in Calabi-Yau manifolds, using moment maps for actions of abelian Lie groups on them. 
The examples include Lagrangian self-shrinkers and translating solitons in the Euclidean spaces. 
Moreover, our method can be applied to construct examples of Lagrangian mean curvature flows in non-flat Calabi-Yau manifolds.  
In particular, we describe Lagrangian mean curvature flows in 4-dimensional Ricci-flat ALE spaces in detail and investigate their singularities.  
\end{abstract}
\section{Introduction}
Special Lagrangian submanifolds in Calabi-Yau manifolds were introduced as special classes of minimal submanifolds in \cite{HL}. 
Since their importance was pointed out in the context of mirror symmetry in \cite{SYZ}, they have been intensively studied. 
As a fundamental tool for constructing special Lagrangian submanifolds, Lagrangian mean curvature flows have been expected. 
In \cite{TY} Thomas-Yau proposed a conjecture, which is reformulated by Joyce \cite{J} recently, that the long-time existence of Lagrangian mean curvature flows are closely related to the stability of Lagrangian submanifolds. 
This conjecture is a central problem in this area. 
Main difficulties in studying mean curvature flows come from singularities of the flow.
In \cite{H} Huisken showed that, if we rescale the mean curvature flow at a type I singularity, the blow-up limit is a self-shrinker. 
It is also known that a translating soliton appears as a certain blow-up limit at another type of singularities. 
Therefore, it is important to study self-shrinkers or translating solitons, which can be considered as local models of singularities of mean curvature flows. 

In the case of Lagrangian mean curvature flows, there are some examples of self-shrinkers or translating solitons as follows. 
Anciaux \cite{A}, Lee-Wang \cite{LW1,LW2} constructed examples of Lagrangian self-shrinkers and self-expanders. 
Later, Joyce-Lee-Tsui \cite{JLT} constructed another type of Lagrangian self-shrinkers and self-expanders as well as Lagrangian translating solitons. 
Moreover, Castro-Lerma \cite{CL} also constructed another type of Lagrangian translating solitons in $\mathbb{C}^2$.
Recently, in \cite{Y}, Yamamoto pointed out that some of these constructions can be interpreted in terms of moment maps for Hamiltonian torus actions on toric Calabi-Yau manifolds and he also constructed examples of generalized Lagrangian mean curvature flows in the sense of \cite{Be} in almost Calabi-Yau manifolds.

In this paper, we construct of various examples of Lagrangian mean curvature flows, using moment maps for actions of abelian Lie groups on Calabi-Yau manifolds. 
Since some parts of Yamamoto's construction rely on toric geometry, we get rid of it so that our construction can be applied to more general cases. 
Roughly speaking, we prove the following. 
See Theorem \ref{construction} for the precise statement. 
\begin{theorem}  \label{construction0}
Let $M$ be a Calabi-Yau manifold and $L$ a special Lagrangian submanifold. 
Suppose there is an isometric and Hamiltonian action of an abelian Lie group $H$   with a moment map $\mu_H \colon M \to \mathfrak{h}^*$. 
Suppose also that the $H$-orbit through $p \in L$ intersects with $L$ orthogonally for each $p \in L$. 
Fix $c_0 \in \mathfrak{h}^*$ and set $V_c=\mu_H^{-1}(c) \cap L$ for $c \in \mathfrak{h}^*$.

Then there exist $a_H \in \mathfrak{h}^*$ and a vector field $\chi$, which is written in terms of $a_H$ explicitly, on $L$ such that the following holds:
\\
$(1)$ Let $\gamma_p \colon [0,T) \to L$ be the integral curve of the vector field $\chi$ with $\gamma_p(0)=p$.
If $p \in V_{c_0}$, then $\gamma_p(t) \in V_{c_t}$ holds, where $c_t= c_0 - t a_H \in \mathfrak{h}^*$. 
\\
$(2)$ Define a family of maps $\{ F_t \colon V_{c_0} \times H \to M \}_{t \in [0,T)}$ by $F_t(p,h)=\gamma_p(t)h$.
If it exists,  then it is a Lagrangian mean curvature flow.
\end{theorem}

Using the above theorem, we construct various examples of Lagrangian mean curvature flows. 
Our basic examples are Lagrangian self-shrinkers, which were already found in \cite{LW2} as well as Lagrangian translating solitons in the Euclidean spaces, which are higher dimensional generalizations of ones given in \cite{CL}. 

Moreover, our construction can be applied to construct Lagrangian mean curvature flows in non-flat Calabi-Yau manifolds.  
4-dimensional Ricci-flat ALE spaces of $A_n$-type are hyperK\"ahler manifolds, on which $S^1$ acts, preserving their hyperK\"ahler structures. 
These $S^1$-actions are extended to $T^2$-actions, which preserve only their K\"ahler structures. 
We use these actions not only for constructing Lagrangian mean curvature flows in these spaces, but also for investigating singularities of the flows. 
As a result, we see that the singularities are type I and determine the blow-up limits at the singularities. 

The condition in Theorem \ref{construction0}, that the $H$-orbit through $p \in L$ intersects with $L$ orthogonally for each $p \in L$, has already appeared in Chapter 9 in \cite{J2}, where Joyce constructed examples of special Lagrangian submanifolds of $\mathbb{C}^m$ by using momant maps. 
His construction is the case $a_H=0$ in Theorem \ref{construction0}. 
Thus our theorem can be considered as a generalization of Joyce's construction of special Lagrangian submanifolds to Lagrangian mean curvature flows. 
The author would like to thank Dominic Joyce for informing him of it. 

Contents of this paper is as follows. 
In Section 2, we state our general construction and prove it. 
In Section 3, we apply the construction to the Euclidean spaces to get Lagrangian self-shrinkers, self-expanders, and translating solitons. 
In Section 4, we apply our construction to hyperK\"ahler manifolds which admit isometric, $I_1$-antiholomorphic and $I_2$-holomorphic involutions. 
As examples of this construction, we describe Lagrangian mean curvature flows in 4-dimensional Ricci-flat ALE spaces of $A_n$-type in Section 5. 
\section{Constructions of Lagrangian mean curvature flows}\label{BC}

Let $(M,I)$ be a complex manifold. 
Here, $I$ is an integrable almost complex structure, which is an automorphism of the tangent bundle $TM$ satisfying $I^2=-1$ and the integrability condition. 
A K\"ahler form $\omega$ on $(M,I)$ is an $I$-invariant symplectic form which induces a Riemannian metric $g(u,v)=\omega (u,I_p v)$ for $p \in M$ and $u,v \in T_pM$. 
\begin{definition}\label{def-CY} 
{\rm A} Calabi-Yau $n$-fold {\rm is a quadruple $(M,I,\omega, \Omega)$ such that $(M,I)$ is an $n$-dimensional complex manifold equipped with a K\"ahler form $\omega$ and a holomorphic volume form $\Omega$ which satisfy the following relation:}
\begin{align}\label{CY}
\frac{\omega^n}{n!}=(-1)^{\frac{n(n-1)}{2}} (\frac{\sqrt{-1}}{2})^n \Omega \wedge \overline{\Omega}.
\end{align}
\end{definition}

It is well known that Calabi-Yau manifolds are Ricci-flat. 
Moreover, if $L$ is an oriented Lagrangian submanifold of a Calabi-Yau manifold $(M,I,\omega, \Omega)$ with the K\"ahler metric $g$,  then there exists a function $\theta \colon L \to \mathbb{R}/2\pi \mathbb{Z}$, which is called the {\it Lagrangian angle}, satisfying
\begin{align*}
\iota^* \Omega = e^{\sqrt{-1}\theta}\mathrm{vol}_{\iota^* g},
\end{align*}
where $\iota \colon L \to M$ is the embedding, and $\mathrm{vol}_{\iota^* g}$ is the volume form of $L$ with respect to the induced metric $\iota^* g$.
Moreover, the mean curvature vector $\mathcal{H}_p$ at $p\in L$ is given by 
\begin{align*}
\mathcal{H}_p = I_{\iota(p)} (\iota_{*p} (\mathrm{grad}_{\iota^* g}\theta)_p ) \in (T_{\iota(p)} L)^\perp,
\end{align*}
where $\mathrm{grad}_{\iota^* g}\theta$ is the gradient of the function $\theta$ with respect to the induced metric $\iota^* g$. 
Then, special Lagrangian submanifolds are defined as follows. 
\begin{definition}
{\rm Let $(M,I,\omega, \Omega)$ be a Calabi-Yau manifold. 
An oriented Lagrangian submanifold $L$ of $M$ is a}  special Lagrangian submanifold with the phase $\theta_0 \in \mathbb{R}/2\pi \mathbb{Z}$ {\rm if its Lagrangian angle $\theta \colon L \to \mathbb{R}/2\pi \mathbb{Z}$ is the constant function $\theta = \theta_0$.}
\end{definition}

Next, we fix our notation on Hamiltonian group actions. 
In this paper, a Lie group $G$ acts on a manifold $M$ from the right. 
We denote the right translation of $g \in G$ by $R_g \colon M \to M$. 
The exponential map is denoted by $\mathrm{Exp}_G \colon \mathfrak{g} \to G$. 
The Lie group $G$ acts on $\mathfrak{g}^*$ by the coadjoint action. 

Suppose that a Lie group $G$ with the Lie algebra $\mathfrak{g}$ acts on a symplectic manifold $(M, \omega)$. 
Then $G$-equivariant map $\mu_G \colon M \to \mathfrak{g}^*$ is a {\it moment map} if
$$-i(\xi^\#) \omega = d \langle \mu_G(\cdot), \xi \rangle$$
holds for each $\xi \in \mathfrak{g}$, where $\xi^\#$ is the vector field on $M$ generated by $\xi$.

Then we have the following theorem. 
\begin{theorem}\label{construction}
Let $(M,I,\omega, \Omega)$ be a connected Calabi-Yau n-fold and $g$  
the K\"ahler metric. 
Let $L$ be a special Lagrangian submanifold with the phase $\theta_0 \in \mathbb{R}/2\pi \mathbb{Z}$.
Suppose that an abelian Lie group $H$ acts on $(M, I, \omega)$, preserving $I$ and $\omega$, with a moment map $\mu_H \colon M \to \mathfrak{h}^*$. 
Suppose also that $\xi^\#_p \ne 0 \in T_pM$ and that $\xi^\#_p \perp T_pL$ for each $p \in L$, $\xi \in \mathfrak{h} \setminus \{ 0 \}$. 
Then the following holds.
\\
$(1)$ There exists $a_H \in \mathfrak{h}^*$ such that $(R_{\mathrm{Exp}_H \xi})^* \Omega = e^{\sqrt{-1} \langle a_H,\xi \rangle}\Omega$ on $M$ for each $\xi \in \mathfrak{h}$.
\\
$(2)$ Set $V_c=\mu_H^{-1}(c) \cap L$ for $c \in \mathfrak{h}^*$.
Then $V_c$ is an orientable submanifold of M of dimension $n - \dim H$ if $V_c$ is not empty. 
Moreover, $T_p V_c \perp I_p \xi^\#_p $ holds for each $p \in V_c$ and $\xi \in \mathfrak{h}$.
\\
$(3)$ 
The map $\phi_c \colon V_c \times H \to M$, which is defined by $\phi_c(p,h)=ph$, is a Lagrangian immersion. 
Moreover, $V_c \times H$ is has the canonical orientation so that $(\phi_c^* \Omega)_{(p, \mathrm{Exp}_H \xi)} = e^{\sqrt{-1} \theta_c} \mathrm{vol}_{\phi_c^* g}$ holds, where $\theta_c(p, \mathrm{Exp}_H \xi) = \langle a_H, \xi \rangle + \theta_0 - \frac{\pi \dim H}{2}$.
\\
$(4)$ Define a map $\tilde{\xi} \colon L \to \mathfrak{h}$ by 
$g_p((\tilde{\xi}(p))^\#_p, \eta^\#_p)= \langle a_H, \eta \rangle$ for $p \in L$ and $\eta \in \mathfrak{h}$. 
Define the vector field $\chi$ on $L$ by $\chi_p=I_p\tilde{\xi}(p)^\#_p \in T_pL$ for $p \in L$.
Then the mean curvature vector $\mathcal{H}^c \in \Gamma(\phi_c^*TM)$ of the map $\phi_c \colon V_c \times H \to M$ is given by $\mathcal{H}^c_{(p,h)}=R_{h*p}\chi_p \in T_{\phi_c(p,h)}M$.
\\
$(5)$ Fix $c_0 \in \mathfrak{h}^*$. 
Suppose that there exists $T >0$ such that, for each $p \in \mu_H^{-1}(c_0) \cap L$, there exists the integral curve $\gamma_p \colon [0,T) \to L$ of the vector field $\chi$ with $\gamma_p(0)=p$.
Then $\gamma_p(t) \in V_{c_t}$ holds for any $p \in L$ and $t \in [0,T)$, where $c_t= c_0 - t a_H \in \mathfrak{h}^*$.
Moreover, the family of maps $\{ F_t \colon V_{c_0} \times H \to M \}_{t \in [0,T)}$, which is defined by $F_t(p,h)=\gamma_p(t)h$,  is a Lagrangian mean curvature flow.
Namely, $\frac{\partial F_t}{\partial t}(p,h)=\mathcal{H}^{c_t}_{(\gamma_p(t),h)}$ holds for $t \in [0,T)$.
\end{theorem}
{\it Proof.}
$(1)$ Fix any $h \in H$. 
Since $R_h^*\Omega$ is a holomorphic $(n,0)$-form, there exists a holomorphic function $f_h \colon M \to \mathbb{C}$ such that $R_h^*\Omega = f_h \Omega$. 
By the relation (\ref{CY}), we have $|f_h(p)|=1$ for each $p \in M$.
Since $M$ is connected, we see that $f_h$ is constant, that is, $f_h(p)=c_h$ for each $p \in M$. 

Moreover, since $R_{h_1 h_2}^* \Omega=R_{h_2}^*(R_{h_1}^*\Omega)$, we have $c_{h_1h_2}=c_{h_1}c_{h_2}$. So, there exists $a_H \in \mathfrak{h}^*$ such that $c_{\mathrm{Exp}_H \xi}=e^{\sqrt{-1} \langle a_H,\xi \rangle}$ for each $\xi \in \mathfrak{h}$.
\\
$(2)$ It is easy to see that $\mathrm{grad}_g \langle \mu_H(\cdot), \xi \rangle = - I \xi^\#$ on $M$. 
Since $L$ is a Lagrangian submanifold of a K\"ahler manifold $(M,I,\omega)$, the complex structure $I_p \colon T_pM \to T_pM$ induces an isomorphism from $T_pL$ to the orthogonal complement $(T_pL)^\perp$ for each $p \in L$.
Moreover, since $\xi^\#_p \ne 0$ and $\xi^\#_p \perp T_pL$ for $p \in L$ by the assumption, we have 
$(\mathrm{grad}_g \langle \mu_H(\cdot), \xi \rangle)_p \in T_pL  \setminus \{ 0 \}$ 
for $p \in L$ and $\xi \in \mathfrak{h} \setminus \{ 0 \}$. 
Therefore, if we set $\nu= \mu_H|_L \colon L \to \mathfrak{h}^*$, 
we have $(d \langle \nu, \xi \rangle)_p \ne 0 \in T^*_p L$ for $p \in L$ and $\xi \in \mathfrak{h} \setminus \{ 0 \}$. 
Since $\langle (d \nu)_p, \xi \rangle = (d \langle \nu, \xi \rangle)_p \ne 0$ for each $\xi \in \mathfrak{h} \setminus \{ 0 \}$, we see that $(d \nu)_p \colon T_pL \to \mathfrak{h}^*$ is surjective for each $p \in L$. 
That is, every point in $L$ is regular for $\nu$. 
This implies that, for any $c \in \mathfrak{h}^*$, the level set $\nu^{-1}(c)= \mu_H^{-1}(c) \cap L$ is a submanifold of dimension $n- \dim H$ if $\nu^{-1}(c)$ is not empty.

For each $v \in T_p(\mu_H^{-1}(c) \cap L)$ and $\xi \in \mathfrak{h}$, we have 
$$0= \langle (d\mu_H)_p(v), \xi \rangle = - \omega_p(\xi^\#_p, v)= -g_p(I_p\xi^\#_p, v).
$$
Thus we see $T_p(\mu_H^{-1}(c) \cap L) \perp I_p \xi^\#_p$ for $p \in \mu_H^{-1}(c) \cap L$ and $\xi \in \mathfrak{h}$.

Since $L$ is oriented, $T_pL$ is an oriented vector space for each $p \in L$. 
Define an isomorphism 
\begin{equation}\label{ori-isom}
\psi_p \colon T_p(\mu_H^{-1}(c) \cap L) \oplus \mathfrak{h} \to T_pL
\end{equation}
by $\psi_p(v,\xi)=v + I_p\xi^\#_p$. 
If we fix an orientation on $\mathfrak{h}$, then we can define an orientation on $T_p(\mu_H^{-1}(c) \cap L)$ so that $\psi_p$ is orientation preserving. 
Thus we see that $\mu_H^{-1}(c) \cap L$ is orientable. 
\\
$(3)$ Define an orientation of $T_{(p,h)}(V_c \times H) \cong T_pV_c \oplus \mathfrak{h}$ for each $(p,h) \in V_c \times H$ so that the map $\psi_p$ in (\ref{ori-isom}) is orientation preserving. 
Thus $V_c \times H$ has the canonical orientation. 
Moreover, we note that, for each $(p,h) \in V_c \times H$ and $(v,\xi) \in T_pV_c \oplus \mathfrak{h} \cong T_{(p,h)}(V_c \times H)$, the following holds:
\begin{align*}
(\phi_c)_{*(p,h)}(v,\xi) = (R_h)_{*p} (v + \xi^\#_p) \in T_{\phi_c(p,h)}M.
\end{align*}

Firstly, we show that $\phi_c$ is an immersion. 
Suppose that $ (\phi_c)_{*(p,h)}(v,\xi) = 0$. 
Then we have $(R_h)_{*p}  (v + \xi^\#_p ) =0 \in T_{ph}M$.  
Since $(R_h)_{*p} \colon T_p M \to T_{ph}M$ is an isomorphism, we have $v + \xi^\#_p =0 \in T_p M$.  
Since $v \perp \xi^\#_p$, we see that $v= \xi^\#_p= 0 \in T_p M$.  
Thus we see that $(v, \xi) = 0 \in T_pV_c \times \mathfrak{h}$, which implies that $\phi_c \colon V_c \times H \to M$ is an immersion.

Secondly, we show $\phi_c^*\omega=0$.
Since $R_h^* \omega=\omega$, we have
\begin{align*}
(\phi_c^* &\omega)_{(p,h)}((v_1, \xi_1), (v_2, \xi_2))\\
&= \omega_{\phi_c(p,h)}((\phi_c)_{*(p,h)}(v_1, \xi_1), (\phi_c)_{*(p,h)}(v_2, \xi_2)) \\
&=(R_h^* \omega)_p(v_1 + (\xi_1)^\#_p, v_2 + (\xi_2)^\#_p) \\
&= \omega_p(v_1,v_2)+\omega_p((\xi_1)^\#_p,(\xi_2)^\#_p) 
+\omega_p(v_1,(\xi_2)^\#_p)+\omega_p((\xi_1)^\#_p,v_2).
\end{align*}
The first term is zero, because $v_i \in T_p L$ for $i=1,2$ and $L$ is a Lagrangian submanifold. 
The second term is zero, because $H$ is abelian. 
For the third term, we have
$$
\omega_p(v_1,(\xi_2)^\#_p) = (d \langle \mu_H(\cdot), \xi_2 \rangle)_p(v_1)=0,
$$
because $v_1 \in T_pV_c$ and $\mu_H$ is constant on $V_c$. 
Similarly, the fourth term is zero. 
Thus we have $\phi_c^* \omega =0$.

Thirdly, we compute $\phi_c^*\Omega$.
Fix an arbitrary $(p,h) \in V_c \times H$ and take $\xi_h \in \mathfrak{h}$ such that $h = \mathrm{Exp}_H \xi_h$.
Since $(\phi_c)_{*(p,h)}(v,0) = (R_h)_{*p}v \in T_{ph}M$ is perpendicular to $(\phi_c)_{*(p,h)}(0,\xi)= (R_h)_{*p}\xi^\#_p \in T_{ph}M$ for each $v \in T_pV_c$ and $\xi \in \mathfrak{h}$, there exists $v_1, \dots, v_{n-m} \in T_pV_c$ and $\xi_1, \dots, \xi_m \in \mathfrak{h}$, where  $m= \dim H$, such that the collection 
$$(v_1, 0), \dots, (v_{n-m},0), (0, \xi_1), \dots, (0,\xi_m) \in T_{(p,h)}(V_c \times H)$$ 
is an oriented orthonormal basis of $T_{(p,h)}(V_c \times H)$ with respect to $(\phi_c^*g)_{(p,h)}$. 
Then, the collection 
$$v_1, \dots, v_{n-m}, I_p(\xi_1)^\#_p, \dots,I_p (\xi_m)^\#_p$$
is an oriented orthonormal basis of $T_p L$.

Recall that $R_h^*\Omega=e^{\sqrt{-1}\langle a_H, \xi_h \rangle} \Omega$ and $\iota_L^* \Omega =e^{\sqrt{-1}\theta_0} \mathrm{vol}_{\iota_L^* g}$, where $\iota_L \colon L \to M$ is the embedding. 
Moreover, since $\xi^\# + \sqrt{-1} I \xi^\#$ is a (0,1)-vector field and $\Omega$ is an $(n,0)$-form on $(M,I)$, we have
$$i(\xi^\#)\Omega=-\sqrt{-1} i( I \xi^\#)\Omega.$$ 
So we have
\begin{align*}
&(\phi_c^* \Omega)_{(p,h)}((v_1, 0), \dots, (v_{n-m},0), (0, \xi_1), \dots, (0,\xi_m))\\
&=\Omega_{ph}((R_h)_{*p} v_1, \dots, (R_h)_{*p} v_{n-m}, (R_h)_{*p}(\xi_1)^\#_p, \dots, (R_h)_{*p}(\xi_m)^\#_p) \\
&=(R_h^*\Omega)_p(v_1, \dots, v_{n-m}, (\xi_1)^\#_p, \dots, (\xi_m)^\#_p) \\
&=e^{\sqrt{-1}\langle a_H, \xi_h \rangle} \Omega_p(v_1, \dots, v_{n-m}, (\xi_1)^\#_p, \dots, (\xi_m)^\#_p) \\
&=e^{\sqrt{-1}\langle a_H, \xi_h \rangle}(-\sqrt{-1})^m \Omega_p(v_1, \dots,v_{n-m}, I_p(\xi_1)^\#_p, \dots,I_p (\xi_m)^\#_p).
\end{align*}
Since the collection $v_1, \dots,v_{n-m}, I_p(\xi_1)^\#_p, \dots,I_p (\xi_m)^\#_p$ is an oriented orthonormal basis of $T_pL$, we have 
\begin{align*}
&(\phi_c^* \Omega)_{(p,h)}((v_1, 0), \dots, (v_{n-m},0), (0, \xi_1), \dots, (0,\xi_m)) \\
& \hspace*{60mm}= e^{\sqrt{-1}\langle a_H, \xi_h \rangle}(-\sqrt{-1})^m e^{\sqrt{-1}\theta_0}.
\end{align*}
Since the collection $(v_1, 0), \dots, (v_{n-m},0), (0, \xi_1), \dots, (0,\xi_m)$ is an oriented orthonormal basis of $T_{(p,h)}(V_c \times \mathfrak{h})$, then we have 
\begin{align*}
\phi_c^* \Omega = e^{\sqrt{-1}\langle a_H, \xi_h \rangle}(-\sqrt{-1})^m e^{\sqrt{-1}\theta_0} \mathrm{vol}_{\phi_c^* g}.
\end{align*}
$(4)$ Let $\mathrm{grad}_{\phi_c^*g}\theta_c$ be the gradient of $\theta_c$ with respect to the metric $\phi_c^*g$. 
Firstly, we show that $(\mathrm{grad}_{\phi_c^*g}\theta_c)_{(p,h)}=(0,\tilde{\xi}(p)) \in T_pV_c \times \mathfrak{h} \cong T_{(p,h)}(V_c \times H)$. 
For each $(v,\eta) \in T_pV_c \times \mathfrak{h} \cong T_{(p,h)}(V_c \times H)$ we have
\begin{align*}
(\phi_c^*g)_{(p,h)}((\mathrm{grad}_{\phi_c^*g}\theta_c)_{(p,h)}, (v,\eta)) 
&=(d \theta_c)_{(p,h)}(v,\eta) \\
&=\langle a_H, \eta \rangle \\
&=g_p((\tilde{\xi}(p))^\#_p, \eta^\#_p) \\
&=g_{ph}((R_h)_{*p}\tilde{\xi}(p)^\#_p, (R_h)_{*p}(v+\eta^\#_p))\\
&=g_{\phi_c (p,h)}((\phi_c)_{*(p,h)}(0,\tilde{\xi}(p)), (\phi_c)_{*(p,h)}(v,\eta))\\
&=(\phi_c^* g)_p((0,\tilde{\xi}(p)), (v,\eta)).
\end{align*}
Thus we see $(\mathrm{grad}_{\phi_c^*g}\theta_c)_{(p,h)}=(0,\tilde{\xi}(p))$. 
Therefore the mean curvature vector $\mathcal{H}^c$ is computed as follows.
\begin{align*}
\mathcal{H}^c_{(p,h)} 
&=I_{\phi_c(p,h)}(\phi^c)_{*(p,h)}(\mathrm{grad}_{\phi_c^*g}\theta_c)_{(p,h)} \\
&=I_{ph}(\phi^c)_{*(p,h)}(0,\tilde{\xi}(p)) \\
&=I_{ph}(R_h)_{*p}(\tilde{\xi}(p)^\#_p) \\
&=(R_h)_{*p}I_{p}(\tilde{\xi}(p)^\#_p) \\
&=(R_h)_{*p}\chi_p.
\end{align*}
$(5)$ Firstly, we show $\gamma_p(t) \in V_{c_t}=\mu_H^{-1}(c_t) \cap L$. 
Since $\gamma_p(t) \in L$ is obvious, it is enough to show $\gamma_p(t) \in \mu_H^{-1}(c_t)$. 
Let $\eta_1, \dots, \eta_m$ be a basis of $\mathfrak{h}$ and $\eta^1, \dots, \eta^m \in \mathfrak{h}^*$ the dual basis. 
Since $\frac{d\gamma_p}{dt}(t)=\chi_{\gamma_p(t)}= I_{\gamma_p(t)}(\tilde{\xi}(\gamma_p(t))^\#_{\gamma_p(t)})$, we have
\begin{align*}
\frac{d}{dt} \mu_H(\gamma_p(t))
&=\frac{d}{dt} \sum_{i=1}^m \langle \mu_H(\gamma_p(t)),\eta_i \rangle \eta^i \\
&=- \sum_{i=1}^m \omega_{\gamma_p(t)}( (\eta_i)^\#_{\gamma_p(t)},\frac{d\gamma_p}{dt}(t)) \eta^i \\
&=- \sum_{i=1}^m \omega_{\gamma_p(t)}( (\eta_i)^\#_{\gamma_p(t)}, I_{\gamma_p(t)}(\tilde{\xi}(\gamma_p(t))^\#_{\gamma_p(t)})) \eta^i \\
&=- \sum_{i=1}^m g_{\gamma_p(t)}( (\eta_i)^\#_{\gamma_p(t)}, \tilde{\xi}(\gamma_p(t))^\#_{\gamma_p(t)}) \eta^i \\
&=- \sum_{i=1}^m \langle a_H,  \eta_i \rangle \eta^i =-a_H.
\end{align*}
Therefore we have $\mu_H(\gamma_p(t))=\mu_H(\gamma_p(0))-ta_H=c_t$.
Thus we see that the image of the map $F_t \colon V_{c_0} \times H \to M$ is contained in the image of the map $\phi_{c_t} \colon V_{c_t} \times H \to M$. 
Since the both maps are immersions, the mean curvature vector of the immersion $F_t$ at $(p,h) \in V_{c_0} \times H$ is $\mathcal{H}^{c_t}_{(\gamma_p(t),h)} \in T_{\gamma_p(t)h}M$. 

On the other hand, we have 
$$
\frac{\partial F_t}{\partial t}(p,h)
=\frac{\partial}{\partial t}(\gamma_p(t)h)
=(R_h)_{* \gamma_p(t)}\frac{d \gamma_p}{dt}(t)
=(R_h)_{* \gamma_p(t)}\chi_{\gamma_p(t)}
=\mathcal{H}^{c_t}_{(\gamma_p(t),h)}.
$$
Thus we finish the proof of Theorem \ref{construction}. \hfill $\Box$

\section{Basic Examples}
In this section we apply Theorem \ref{construction} and give basic examples of Lagrangian mean curvature flows in the Euclidean spaces.
\subsection{Self-shrinkers and self-expanders}\label{3-1}

Let $M=\mathbb{C}^d$ equipped with the standard complex structure $I$, the standard K\"ahler form $\omega$ and the standard holomorphic volume form $\Omega$. 
The K\"ahler metric is denoted by $g$. 

Let $H$ be $S^1= \{ \zeta \in \mathbb{C}~|~|\zeta|=1\}$ and $\mathfrak{h}$ its Lie algebra. 
Take $\xi_0 \in \mathfrak{h}$ such that  $\mathrm{Exp}_H t \xi_0 =e^{\sqrt{-1}t} \in H$. 
Take $\xi^0 \in \mathfrak{h}^*$ such that $\langle \xi^0, \xi_0 \rangle =1$.
The action of $H$ on $M$ is defined by 
$$
z\mathrm{Exp}_H t \xi_0 =(z_1e^{\sqrt{-1}\lambda_1 t}, \dots, z_d e^{\sqrt{-1}\lambda_d t}),
$$
where $z=(z_1, \dots, z_d) \in M$ and $\lambda_i \in \mathbb{Z} \setminus \{ 0 \}$ for $i= 1, \dots, d$. 
The moment map $\mu_H \colon M \to \mathfrak{h}^*$ is given by $\mu_H(z)= \frac{1}{2} (\sum_{i=1}^d \lambda_i |z_i|^2)\xi^0$. 
We also have 
$$
R_{\mathrm{Exp}_H \xi}^* \Omega=e^{\sqrt{-1} \langle a_H, \xi \rangle}\Omega, \quad \text{where $\xi \in \mathfrak{h}$ and $a_H=(\lambda_1 + \dots + \lambda_d)\xi^0 \in \mathfrak{h}^*$.}
$$

Note that $(\xi_0)^\#_x=(\sqrt{-1}\lambda_1 x_1, \dots, \sqrt{-1}\lambda_d x_d) \in \mathbb{C}^d \cong T_xM$ for each $x =(x_1, \dots, x_d) \in \mathbb{R}^d$. 
So, if we set $L= \mathbb{R}^d \setminus \{ (0, \dots, 0) \}$, then $L$ is a special Lagrangian submanifold of $(M,I,\omega,\Omega)$, and $(\xi_0)^\#_x \ne 0 \in T_xM$ and $(\xi_0)^\#_x \perp T_xL$ hold for each $x \in L$. 

\begin{lemma}\label{chi1}
Under the above setting, the vector field $\chi$ on L in Theorem \ref{construction} $(4)$ is given by 
$$
\chi_x=\frac{- \sum_{i=1}^d \lambda_i}{\sum_{i=1}^d \lambda_i^2 x_i^2}(\lambda_1 x_1, \dots, \lambda_d x_d) \in \mathbb{R}^d \cong T_xL, 
$$
where $x=(x_1, \dots, x_d) \in L$. 
In particular, $\lim_{|x| \to \infty}|\chi_x| = 0$ holds.
\end{lemma}
{\it Proof.} Note that $(\xi_0)^\#_x=(\sqrt{-1}\lambda_1 x_1, \dots, \sqrt{-1}\lambda_d x_d)$ for each $x \in L$ and that $\tilde{\xi}(x) \in \mathfrak{h}$ is defined by 
$g_x((\tilde{\xi}(x))^\#_x, (\xi_0)^\#_x)= \langle a_H, \xi_0 \rangle$.
If we set $\tilde{\xi}(x)=\alpha_x \xi_0$, where $\alpha_x \in \mathbb{R}$, then we have 
$g_x(\alpha_x(\xi_0)^\#_x,(\xi_0)^\#_x)= \langle a_H, \xi_0 \rangle$. 
So we have 
$$
\alpha_x = \frac{\langle a_H, \xi_0 \rangle}{g_x((\xi_0)^\#_x,(\xi_0)^\#_x)}= \frac{\sum_{i=1}^d \lambda_i}{\sum_{i=1}^d \lambda_i^2 x_i^2}.
$$
Since $\chi_x= I_x (\tilde{\xi}(x))^\#_x = \sqrt{-1} \alpha_x (\xi_0)^\#_x$, we finish the proof.  \hfill $\Box$

Then we have the following.

\begin{proposition}\label{prop-existence}
Fix $c_0 \in \mathfrak{h}^* \setminus \{ 0 \}$ and set $c_t=c_0 -ta_H$. 
Suppose $c_t \ne 0$ for $t \in [0,T)$.
Then the family of maps $\{ F_t \colon (\mu_H^{-1}(c_0) \cap L) \times H \to M \}_{t \in [0,T)}$, which is defined by $F_t(p,h)=\gamma_p(t)h$, is a Lagrangian mean curvature flow.
\end{proposition}
{\it Proof.} First we show that, for each $p \in \mu_H^{-1}(c_0) \cap L$,  the integral curve $\{ \gamma_p(t) \}$ of  the vector field $\chi$ with $\gamma_p(0)=p$ exists for $t \in [0,T)$. 
Assume that, for some $p_0 \in \mu_H^{-1}(c_0) \cap L$, $\gamma_{p_0}(t)$ exists only for $t \in [0,t_0)$, where $0<t_0 < T$.
By Lemma \ref{chi1}, we see that $| \frac{d \gamma_{p_0}}{dt}(t)|$ is uniformly bounded for $t \in [0,t_0)$. 
Since $c_{t_0} \ne 0$, $\mu_H^{-1}(c_{t_0}) \cap L$ is closed in $\mathbb{R}^n$. 
So it is easy to see that $\lim_{t \to t_0}\gamma_{p_0}(t)$ exists in $\mu_H^{-1}(c_{t_0}) \cap L$ and that $\gamma_{p_0}(t)$ extends to $t \in [0,t_0 + \epsilon)$ for some $\epsilon > 0$. 
This contradicts to the assumption. 
Thus we see that, for each $p \in \mu_H^{-1}(c_0) \cap L$,  the integral curve $\{ \gamma_p(t) \}$ of  the vector field $\chi$ with $\gamma_p(0)=p$ exists for $t \in [0,T)$. 

Then the claim follows from Theorem \ref{construction} $(5)$.  \hfill $\Box$

The above examples are already known. 
In \cite{LW2} Lee and Wang found them as examples of Lagrangian self-shrinkers or self-expanders. 
Here we reprove their results in our context.

\begin{proposition}{\rm\cite{LW2}}
Fix $c \in \mathfrak{h}^* \setminus \{ 0 \}$ and suppose $\mu_H^{-1}(c) \cap L \ne \emptyset$. 
Define a map $\phi_c \colon (\mu_H^{-1}(c)\cap L) \times H \to M$ by $\phi_c(x,h)=xh$.
Let $\phi_c(x,h)^\perp$ be the normal component of the position vector $\phi_c(x,h)$ in $T_{\phi_c(x,h)}M$. 
Then the mean curvature vector $\mathcal{H}^c \in \Gamma(\phi_c^*TM)$ is given by 
$$
\mathcal{H}^c_{(x,h)}= \alpha_c \phi_c(x,h)^\perp \in T_{\phi_c(x,h)}M, \quad \text{where $\alpha_c= \frac{-\sum_{i=1}^d \lambda_i}{2 \langle c,\xi_0 \rangle}$.}
$$
That is, $\phi_c$ is a self-shrinker if $\alpha_c<0$, and a self-expander if $\alpha_c >0$.
\end{proposition}
{\it Proof.}  
Recall $\langle \mu_H(x), \xi_0 \rangle = \frac{1}{2} \sum_{i=1}^d \lambda_i x_i^2$ for  $x=(x_1, \dots, x_d) \in L$.
Since the gradient of $\langle \mu_H( \cdot ), \xi_0 \rangle|_L$ is a normal vector field of $\mu_H^{-1}(c) \cap L$ in $L$,  the unit normal vector filed $\nu$ of $\mu_H^{-1}(c)\cap L$ in $L$ is given by
$$
\nu_x=\frac{1}{\sqrt{\sum_{i=1}^d \lambda_i^2 x_i^2}}(\lambda_1 x_1, \dots, \lambda_d x_d).
$$
For $x \in \mu_H^{-1}(c) \cap L$, fix an orthonormal basis $v_1, \dots, v_d$ of $(T_x \mathrm{Im}\phi_c)^\perp$ satisfying $v_1=\nu_x$.
Then we have $v_i \perp T_x L$, in particular $g_x(x,v_i)=0$, for $i=2, \dots, d$.
Moreover, $(R_h)_{*x}v_1, \dots, (R_h)_{*x}v_d$ is an orthonormal basis of $(T_{\phi_c(x,h)} \mathrm{Im}\phi_c)^\perp$.
So we have, for $(x,h) \in (\mu_H^{-1}(c)\cap L) \times H$,
\begin{align*}
\phi_c(x,h)^\perp
&=\sum_{i=1}^d g_{xh}(\phi_c(x,h), (R_h)_{*x}v_i)(R_h)_{*x}v_i \\
&=\sum_{i=1}^d g_{x}(x, v_i)(R_h)_{*x}v_i \\
&=(R_h)_{*x}(\sum_{i=1}^d g_{x}(x, v_i)v_i) \\
&=(R_h)_{*x}(g_{x}(x, v_1)v_1) \\
&=(R_h)_{*x}(\frac{\sum_{i=1}^d \lambda_i x_i^2}{\sqrt{\sum_{i=1}^d \lambda_i^2 x_i^2}}\frac{1}{\sqrt{\sum_{i=1}^d \lambda_i^2 x_i^2}}(\lambda_1 x_1, \dots, \lambda_d x_d)) \\
&=(R_h)_{*x}(\frac{2\langle c, \xi_0 \rangle}{\sum_{i=1}^d \lambda_i^2 x_i^2}(\lambda_1 x_1, \dots, \lambda_d x_d)) .
\end{align*}
Together with Lemma \ref{chi1} and Theorem \ref{construction} $(4)$,  we have $$\alpha_c \phi_c(x,h)^\perp =(R_h)_{*x}\chi_x=\mathcal{H}^c_{(x,h)}.$$
\hfill $\Box$

Set $L_t = \{ xh ~|~ x \in \mu_H^{-1}(c_t) \cap \mathbb{R}^d, h \in H\}$ for $t \in \mathbb{R}$, where $c_t=c_0 -ta_H \in \mathfrak{h}^*$ as above.  
Note that $L_t= \mathrm{Im} \phi_{c_t}$ if $c_t \ne 0$. 
Lee and Wang also proved that $\{ L_t \}_{t \in \mathbb{R}}$ forms an eternal solution for Brakke flow in \cite{LW1,LW2}.

\subsection{Translating solitons}

Let  $M=\mathbb{C}^{d+1}$ equipped with the standard complex structure $I$, the standard K\"ahler form $\omega$ and the standard holomorphic volume form $\Omega$. 
The K\"ahler metric is denoted by $g$. 

Let $H=\mathbb{R}$ and $\mathfrak{h}$ its Lie algebra. 
Take $\xi_0 \in \mathfrak{h}$ such that  $\mathrm{Exp}_H t \xi_0 =t \in H$. 
Take $\xi^0 \in \mathfrak{h}^*$ such that $\langle \xi^0, \xi_0 \rangle =1$.
The action of $H$ on $M$ is defined by 
$$
z\mathrm{Exp}_H t \xi_0 =(z_1e^{\sqrt{-1}\lambda_1 t}, \dots, z_d e^{\sqrt{-1}\lambda_d t}, z_{d+1}+ \sqrt{-1}t),
$$
where $z=(z_1, \dots, z_{d+1}) \in M$ and $\lambda_i \in \mathbb{R}$ for $i= 1, \dots, d$. 
The moment map $\mu_H \colon M \to \mathfrak{h}^*$ is given by $\mu_H(z)= (\frac{1}{2} \sum_{i=1}^d \lambda_i |z_i|^2 + \mathrm{Re} z_{d+1})\xi^0$, where $\mathrm{Re} z_{d+1}$ is the real part of $z_{d+1} \in \mathbb{C}$. 
We also have 
$$
R_{\mathrm{Exp}_H \xi}^* \Omega=e^{\sqrt{-1} \langle a_H, \xi \rangle}\Omega, \quad \text{where $\xi \in \mathfrak{h}$ and $a_H=(\lambda_1 + \dots + \lambda_d)\xi^0 \in \mathfrak{h}^*$.}
$$

Set $L= \mathbb{R}^{d+1}$. 
Then $L$ is a special Lagrangian submanifold of $(M, I,\omega,\Omega)$.
Since $(\xi_0)^\#_x=(\sqrt{-1}\lambda_1 x_1, \dots, \sqrt{-1}\lambda_d x_d, \sqrt{-1}) \in \mathbb{C}^{d+1} \cong T_xM$ for each $x \in L$, we have $(\xi_0)^\#_x \ne 0 \in T_xM$ and $(\xi_0)^\#_x \perp T_xL$. 
\begin{lemma}\label{chi2}
Under the above setting, the vector field $\chi$ on L in Theorem \ref{construction} $(4)$ is given by 
$$
\chi_x=\frac{- \sum_{i=1}^d \lambda_i}{1+\sum_{i=1}^d \lambda_i^2 x_i^2}(\lambda_1 x_1, \dots, \lambda_d x_d,1) \in \mathbb{R}^{d+1} \cong T_xL, 
$$
where $x=(x_1, \dots, x_{d+1}) \in L$. 
In particular, $|\chi_x|$ is uniformly bounded on $L$.
\end{lemma}
{\it Proof.} Note that $(\xi_0)^\#_x=(\sqrt{-1}\lambda_1 x_1, \dots, \sqrt{-1}\lambda_d x_d, \sqrt{-1})$ for each $x \in L$ and that $\tilde{\xi}(x) \in \mathfrak{h}$ is defined by 
$g_x((\tilde{\xi}(x))^\#_x, (\xi_0)^\#_x)= \langle a_H, \xi_0 \rangle$.
If we set $\tilde{\xi}(x)=\alpha_x \xi_0$, where $\alpha_x \in \mathbb{R}$, then we have 
$g_x(\alpha_x(\xi_0)^\#_x,(\xi_0)^\#_x)= \langle a_H, \xi_0 \rangle$. 
So we have 
$$
\alpha_x = \frac{\langle a_H, \xi_0 \rangle}{g_x((\xi_0)^\#_x,(\xi_0)^\#_x)}= \frac{\sum_{i=1}^d \lambda_i}{1+\sum_{i=1}^d \lambda_i^2 x_i^2}.
$$
Since $\chi_x= I_x (\tilde{\xi}(x))^\#_x = \sqrt{-1} \alpha_x (\xi_0)^\#_x$, we finish the proof.  \hfill $\Box$

In contrast to the case of Proposition 1, $\mu_H^{-1}(c_t) \cap L$ is closed in $\mathbb{R}^n$ for any $t \in \mathbb{R}$ in this case. 
So, due to the above lemma, by the same argument as in the proof of Proposition \ref{prop-existence}, we see that, for any $p \in L$, the integral curve $\{ \gamma_p(t) \}$ of  the vector field $\chi$ with $\gamma_p(0)=p$ exists for $t \in \mathbb{R}$ . 
Theorem \ref{construction} implies the following.

\begin{proposition}
Fix $c_0 \in \mathfrak{h}^*$ and set $c_t=c_0 -ta_H$ for $t \in \mathbb{R}$. 
Then the family of maps $\{ F_t \colon (\mu_H^{-1}(c_0) \cap L) \times H \to M \}_{t \in \mathbb{R}}$, which is defined by $F_t(p,h)=\gamma_p(t)h$, is a Lagrangian mean curvature flow.
\end{proposition}

The next proposition shows that the above examples are Lagrangian translating solitons.  
These are higher dimensional generalization of the examples given in \cite{CL}.

\begin{proposition}
Fix $c \in \mathfrak{h}^*$. 
Define a map $\phi_c \colon (\mu_H^{-1}(c)\cap L) \times H \to M$ by $\phi_c(x,h)=xh$. 
Set $u=(0, \dots, 0,-\sum_{i=1}^d \lambda_i) \in \mathbb{R}^{d+1}$.
Let $u^\perp_{(x,h)}$ be the normal component of the vector $u$ in $T_{\phi_c(x,h)}M$. 
Then the mean curvature vector is given by $\mathcal{H}^c_{(x,h)}= u^\perp_{(x,h)} \in T_{\phi_c(x,h)}M$.
That is, $\phi_c$ is a translating soliton.
\end{proposition}
{\it Proof.}  
Recall $\langle \mu_H(x), \xi_0 \rangle = \frac{1}{2} \sum_{i=1}^d \lambda_i x_i^2 + x_{d+1}$ for  $x=(x_1, \dots, x_{d+1}) \in L$.
Since the gradient of $\langle \mu_H( \cdot ), \xi_0 \rangle|_L$ is a normal vector field of $\mu_H^{-1}(c) \cap L$ in $L$,  the unit normal vector filed $\nu$ of $\mu_H^{-1}(c) \cap L$ in $L$ is given by
$$
\nu_x=\frac{1}{\sqrt{1+\sum_{i=1}^d \lambda_i^2 x_i^2}}(\lambda_1 x_1, \dots, \lambda_d x_d,1).
$$
For $x \in \mu_H^{-1}(c) \cap L$, fix an orthonormal basis $v_1, \dots, v_{d+1}$ of $(T_x \mathrm{Im}\phi_c)^\perp$ satisfying $v_1=\nu_x$.
Then we have $v_i \perp T_x L$, in particular $g_x(u, v_i)=0$, for $i=2, \dots, d+1$.
Moreover, $(R_h)_{*x}v_1, \dots, (R_h)_{*x}v_{d+1}$ is an orthonormal basis of $(T_{\phi_c(x,h)} \mathrm{Im}\phi_c)^\perp$.
If we note $(R_h)_{*x}u=u$,  we have, for $(x,h) \in (\mu_H^{-1}(c)\cap L) \times H$,
\begin{align*}
u^\perp_{(x,h)}
&=\sum_{i=1}^{d+1} g_{xh}(u, (R_h)_{*x}v_i)(R_h)_{*x}v_i \\
&=\sum_{i=1}^{d+1} g_{x}(u, v_i)(R_h)_{*x}v_i \\
&=(R_h)_{*x}(\sum_{i=1}^{d+1} g_{x}(u, v_i)v_i) \\
&=(R_h)_{*x}(g_{x}(u, v_1)v_1) \\
&=(R_h)_{*x}(\frac{-\sum_{i=1}^d \lambda_i}{\sqrt{1+\sum_{i=1}^d \lambda_i^2 x_i^2}}\frac{1}{\sqrt{1+\sum_{i=1}^d \lambda_i^2 x_i^2}}(\lambda_1 x_1, \dots, \lambda_d x_d,1)) \\
&=(R_h)_{*x}\chi_x=\mathcal{H}^c_{(x,h)},
\end{align*}
where the last two equalities follow from Lemma \ref{chi2} and Theorem \ref{construction} $(4)$, respectively.
\hfill $\Box$

\section{Lagrangian mean curvature flows in hyperK\"ahler manifolds}
 In this section we give examples of special Lagrangian submanifolds, which satisfies the conditions in Theorem \ref{construction}, in non-flat Calabi-Yau manifolds. 
\begin{lemma}\label{K-involution}
Let $(M,I,\omega)$ be a K\"ahler manifold.
Let $\sigma \colon M \to M$ be an isometric, anti-holomorphic involution, whose fixed point set is denoted by $M^\sigma$. 
Then the following holds.
\\
$(1)$ $\sigma^* \omega = - \omega.$
\\
$(2)$ $M^\sigma$ is a Lagrangian submanifold. 
\\
$(3)$ For $p \in M^\sigma$, the set of eigenvalues of $\sigma_{*p} \colon T_pM \to T_pM$  is $\{1,-1 \}$. 
The eigenspaces $V(\lambda)$, corresponding to the eigenvalues $\lambda=\pm 1$, are given by $V_p(1)= T_p(M^\sigma)$ and 
$V_p(-1)= T_p(M^\sigma)^\perp$, respectively.
\end{lemma}
{\it Proof.} $(1)$ Let $g$ be the K\"ahler metric on $(M,I,\omega)$.
For $p \in M$ and $v,w \in T_pM$, we have 
\begin{align*}
(\sigma^* \omega)_p&(v,w)=\omega_{\sigma(p)}(\sigma_{*p}v, \sigma_{*p}w)
=g_{\sigma(p)}(I_{\sigma(p)} \sigma_{*p}v, \sigma_{*p}w)
\\ &=-g_{\sigma(p)}(\sigma_{*p}I_p v, \sigma_{*p}w)
=-(\sigma^*g)_p(I_p v, w)=-g_p(I_p v, w)=-\omega_p(v,w).
\end{align*}
$(3)$ By the slice theorem, we see that $M^\sigma$ is a submanifold.  
Since $\sigma_{*p} \colon T_pM \to T_pM$ is an isometric involution for $p \in M^\sigma$, we have an orthogonal decomposition $T_pM = V_p(1) \oplus V_p(-1)$ and $T_p (M^\sigma)=V_p(1)$ for $p \in M^\sigma$. 
Then $(3)$ follows.  
\\
$(2)$ If $p \in M^\sigma$ and $v,w \in V_p(\lambda)$, where $\lambda = \pm 1$, then we have
\begin{align*}
\omega_p(v,w)=-(\sigma^* \omega)_p(v,w)=-\omega_p(\sigma_{*p}v, \sigma_{*p}w)
=-\omega_p(\lambda v, \lambda w)=-\omega_p(v,w).
\end{align*}
So we have $\omega_p(v,w)=0$ for $v,w \in V_p(\lambda)$. 
Therefore we have $\dim V_p(\lambda) \le n$, where $n$ is the complex dimension of $(M,I)$. 
Since $T_pM = V_p(1) \oplus V_p(-1)$, we have $\dim V_p(1) = \dim V_p(-1) = n$. 
Thus we see that $M^\sigma$ is a Lagrangian submanifold. \hfill $\Box$

Let us recall the notion of hyperK\"ahler manifolds.
\begin{definition}
{\em A} hyperK\"ahler manifold {\em is a collection $(M,g, I_1,I_2,I_3)$ such that $(M,g)$ is a $4n$-dimensional Riemannian manifold with three complex structures $I_1,I_2,I_3$ which satisfies the following properties:
\\ $(i)$ $I_1,I_2,I_3$ satisfy the quaternion relation, that is, $I_1 I_2=-I_2 I_1=I_3$.
\\ $(ii)$ $g$ is a K\"ahler metric with respect to each complex structure $I_1,I_2,I_3$.}
\end{definition}
Let $\omega_j$ be the K\"ahler form of $(M,g,I_j)$ for $j=1,2,3$. 
Then $\omega_\mathbb{C}= \omega_2+ \sqrt{-1}\omega_3$ is a holomorphic symplectic form and $\displaystyle \Omega=\frac{\omega_\mathbb{C}^n}{n!}$ is a holomorphic volume form of $(M,I_1)$. 
Moreover, $(M,I_1, \omega_1, \Omega)$ is a Calabi-Yau $2n$-fold in the sense of Definition \ref{def-CY}.
\begin{proposition}\label{construction2}
Let $(M,g,I_1,I_2,I_3)$ be a $4n$-dimensional hyperK\"ahler manifold, admitting an isometric involution $\sigma \colon M \to M$, which is anti-holomorphic with respect to $I_1$, and holomorphic with respect to $I_2$, respectively. 
Denote the set of fixed points of $\sigma$ by $M^\sigma$, and the embedding by $\iota \colon M^\sigma \to M$. 
Let $\omega_j$ be the K\"ahler form for $j=1,2,3$. 
Set $\omega_\mathbb{C}= \omega_2+ \sqrt{-1}\omega_3$ and $\displaystyle \Omega=\frac{\omega_\mathbb{C}^n}{n!}$. 
Then the following holds.
\\
$(1)$ $\sigma^* \omega_1 = - \omega_1$, $\sigma^* \omega_2 = \omega_2$ and $\sigma^* \omega_3 = - \omega_3$.
In particular, $M^\sigma$ is a Lagrangian submanifold of $(M, \omega_j)$ for $j=1,3$. 
\\
$(2)$ $M^\sigma$ is a complex submanifold of $(M, I_2)$.
In particular, $M^\sigma$ is oriented. 
\\
$(3)$ $\iota^*\Omega = \mathrm{vol}_{\iota^* g}$ holds. $M^\sigma$ is a special Lagrangian submanifold of $(M, I_1, \omega_1, \Omega)$.
\\
$(4)$ In addition, suppose that an abelian Lie group $H$ acts on $M$, preserving $I_1$ and $\omega_1$,  with a moment map $\mu_H \colon M \to \mathfrak{h}^*$, and that $\sigma(ph)=\sigma(p)h^{-1}$ holds for each $p \in M$ and $h \in H$.
Set $L= \{ p \in M^\sigma~|~ \xi^\#_p \ne 0 ~\text{for each $\xi \in \mathfrak{h} \setminus \{0 \}$} \}$. 
Then $\xi^\#_p \perp T_pL$ for each $p \in L$, $\xi \in \mathfrak{h} \setminus \{0 \}$ holds. 
That is, $L$ satisfies the conditions in Theorem \ref{construction}.
\end{proposition}
{\it Proof.} 
$(1)$ Note that $\sigma \colon M \to M$ is an anti-holomorphic involution with respect to $I_3=I_1 I_2$. 
The claim follows from Lemma \ref{K-involution}.
\\
$(2)$ Since $(I_2)_p \colon T_pM \to T_pM$ commutes with $\sigma_{*p} \colon T_pM \to T_pM$ for $p \in M^\sigma$, $(I_2)_p$ preserves the eigenspaces $V_p(\lambda)$ of $\sigma_{*p}$ for $\lambda= \pm 1$.
Since $T_p(M^\sigma)=V_p(1)$ by Lemma \ref{K-involution}, the claim follows. 
\\
$(3)$ Since $\iota^*\omega_2$ is the K\"ahler form of $(M^\sigma, I_2)$ by $(2)$, we see that $\frac{(\iota^*\omega_2)^n}{n!}$ is the volume form $\mathrm{vol}_{\iota^* g}$ of $M^\sigma$. 
On the other hand, since $\iota^*\omega_\mathbb{C}= \iota^*\omega_2$ by $(1)$ , we have 
$$\iota^*\Omega=\iota^*(\frac{\omega_\mathbb{C}^n}{n!})=\frac{(\iota^*\omega_\mathbb{C})^n}{n!}=\frac{(\iota^*\omega_2)^n}{n!}=\mathrm{vol}_{\iota^* g}.
$$ 
$(4)$ Since $\sigma(ph)=\sigma(p)h^{-1}$ holds for each $p \in M$ and $h \in H$, we have $\sigma_{*p}\xi^\#_p = - \xi^\#_p$ for each $p \in M^\sigma$ and $\xi \in \mathfrak{h}$. 
By Lemma \ref{K-involution}, we have $\xi^\#_p \perp T_p(M^\sigma)$ for each $p \in M^\sigma$ and $\xi \in \mathfrak{h}$. \hfill $\Box$
\section{Lagrangian mean curvature flows in Ricci-flat ALE spaces of type $A_n$}
In this section, as an application of Proposition \ref{construction2}, we construct Lagrangian mean curvature flows in 4-dimensional Ricci-flat ALE spaces of $A_n$-type and study their properties. 

\subsection{Ricci-flat ALE spaces of type $A_n$}

In this subsection we construct 4-dimensional Ricci-flat ALE spaces of $A_n$-type.
Our construction seems to be slightly different from the original one due to \cite{Kr}, but essentially the same. 

In the following, we identify the quaternionic vector space $\mathbb{H}^{n+1}$ with the product $\mathbb{C}^{n+1}\times \mathbb{C}^{n+1}$. 
The standard hyperK\"ahler structure $(g,I_1,I_2,I_3)$ on $\mathbb{H}^{n+1}$ is defined by
\begin{align*}
&g((z,w),(z^\prime,w^\prime))=\mathrm{Re}(^t\!z\overline{z^\prime}+^t\!w\overline{w^\prime}),\\
&I_1(z,w)=(\sqrt{-1}z,\sqrt{-1}w), I_2(z,w)=(-\overline{w},\overline{z}),I_3(z,w)=(-\sqrt{-1}\overline{w},\sqrt{-1}\overline{z}),
\end{align*}
where $z,z^\prime, w,w^\prime \in \mathbb{C}^{n+1}$ are column vectors. 
The K\"ahler forms $\omega_i$ corresponding to the complex structure $I_i$ for  $i=1,2,3$ are given by
\begin{align*}
&\omega_1((z,w),(z^\prime, w^\prime))=\mathrm{Re}\{ \sqrt{-1}(^t\!z\overline{z^\prime}+ \! ^t\!w\overline{w^\prime})\},\\
&(\omega_2 + \sqrt{-1} \omega_3)((z,w),(z^\prime, w^\prime))= \! ^t\!z w^\prime- \! ^t\!z^\prime w.
\end{align*}
The right action of the torus $T^{n+1} = \{ \zeta=(\zeta_0, \cdots, \zeta_n) \in \mathbb{C}^{n+1}~|~ |\zeta_i|=1\quad \text{for}~~  i= 0, \dots, n \}$
on $\mathbb{H}^{n+1}=\mathbb{C}^{n+1}\times \mathbb{C}^{n+1}$ is defined by
\begin{align}\label{equ-action}
( \begin{pmatrix}z_0 \\ \vdots \\ z_n \end{pmatrix} , 
 \begin{pmatrix}w_0 \\ \vdots \\ w_n \end{pmatrix}) \zeta 
=
( \begin{pmatrix}z_0 \zeta_0 \\ \vdots \\ z_n \zeta_n \end{pmatrix} , 
 \begin{pmatrix}w_0 \zeta_0^{-1} \\ \vdots \\ w_n \zeta_n^{-1} \end{pmatrix}) .
\end{align}
Define a group homomorphism $\rho \colon T^{n+1} \to T^1$ by $\rho(\zeta)=\zeta_0 \dots \zeta_n$. 
Denote  the kernel of $\rho$ by $K$, and let $\iota \colon K \to T^{n+1}$ be the embedding. 
Then we have the following exact sequence of abelian Lie groups:
\begin{align}\label{exact}
1 \longrightarrow K \overset{\iota}{\longrightarrow} T^{n+1} \overset{\rho}{\longrightarrow}T^1 \longrightarrow 1.
\end{align}
We also have the corresponding exact sequences of the Lie algebras and their dual spaces:
\begin{align*}
0 \longrightarrow \mathfrak{k} \overset{\iota_*}{\longrightarrow} \mathfrak{t}^{n+1} \overset{\rho_*}{\longrightarrow}\mathfrak{t}^1 \longrightarrow 0, \quad
0 \longleftarrow \mathfrak{k}^* \overset{\iota^*}{\longleftarrow} (\mathfrak{t}^{n+1})^* \overset{\rho^*}{\longleftarrow}(\mathfrak{t}^1)^* \longleftarrow 0.
\end{align*}
Let $e_0, \dots, e_n$ be the standard basis of $\mathfrak{t}^{n+1}$ and $e^0, \dots, e^n \in  (\mathfrak{t}^{n+1})^*$ the dual basis. 
We also denote the standard basis of $\mathfrak{t}^1$ by $p_1$ and the dual basis by $p^1 \in (\mathfrak{t}^1)^*$. 
Then we have 
\begin{align*}
\rho_*e_0= \cdots = \rho_*e_n= p_1, \quad \rho^* p^1=e^0 + \cdots + e^n .
\end{align*}
If we set $f_i=e_i -e_{i-1}$ for $i=1,\dots, n$, then $f_1, \dots, f_n$ is a basis of $\mathfrak{k}$.
The dual basis $f^1, \dots, f^n \in \mathfrak{k}^*$ satisfies the following;
\begin{align}\label{equ-f}
\iota^*e^0= -f^1, \quad \iota^*e^1 = f^1 -f^2, \quad \cdots, \quad \iota^*e^{n-1}=f^{n-1}-f^n, \quad \iota^*e^n=f^n.
\end{align}
The induced $K$-action on $\mathbb{H}^{n+1}$ admits a hyperK\"ahler moment map $$\mu_K=(\mu_K^1,\mu_K^2, \mu_K^3) \colon \mathbb{H}^{n+1} \to \mathfrak{k}^* \otimes \mathrm{Im}\mathbb{H},$$ which is given by
\begin{align*}
&\mu_K^1(z,w)=\frac{1}{2} \sum_{i=0}^n(|z_i|^2 -|w_i|^2)\iota^*e^i=\frac{1}{2} \sum_{i=1}^n (|z_i|^2 -|w_i|^2-|z_{i-1}|^2 +|w_{i-1}|^2)f^i, \\
& (\mu_K^2 + \sqrt{-1}\mu_K^3)(z,w)=-\sqrt{-1} \sum_{i=0}^n z_i w_i\iota^*e^i= -\sqrt{-1} \sum_{i=1}^n (z_i w_i -z_{i-1} w_{i-1} )f^i.
\end{align*}
Define the codimension one subspaces $W_{i,j}$ of $\mathfrak{k}^*$ for $0 \le i < j \le n$ by 
\begin{align*}
W_{i,j} &= \mathrm{span}\{ \iota^*e^0, \dots, \widehat{\iota^*e^i}, \dots, \widehat{\iota^*e^j},  \dots, \iota^*e^n \} \\
&= \{ \alpha \in \mathfrak{k}^*~|~ \langle \alpha, \sum_{k=i+1}^j f_k \rangle =0 \},
\end{align*}
where the second equality follows from $\displaystyle \sum_{k=i+1}^j f_k =e_j -e_i$.
\begin{proposition}\label{prop-rv} 
The set of regular values of the hyperK\"ahler moment map 
\\
$\mu_K \colon \mathbb{H}^{n+1} \to \mathfrak{k}^* \otimes \mathrm{Im}\mathbb{H}$ is
\begin{align*}
\mathfrak{k}^* \otimes \mathrm{Im}\mathbb{H} \setminus \bigcup_{0 \le i < j \le n} W_{i,j}\otimes \mathrm{Im}\mathbb{H}.
\end{align*}
\end{proposition}
{\it Proof.}
In general, let us define a map $\nu \colon \mathbb{C}^2 \to \mathbb{R} \times \mathbb{C}$ by $\nu(u,v)=(|u|^2 - |v|^2, uv)$. 
Then it is easy to observe that $\nu$ is surjective, and that the set of regular points of $\nu$ is $\mathbb{C}^2 \setminus \{ (0,0)\}$.
 
By this observation, the image of the differential of $\mu_K$ at $(z,w) \in \mathbb{H}^{n+1}$ is 
\begin{align*}
\sum_{j \in J_{(z,w)}}\iota^* e^j \otimes \mathrm{Im}\mathbb{H}, \quad \text{where $J_{(z,w)}= \{ j ~|~ (z_j,w_j) \ne (0,0) \}$.}
\end{align*}
Therefore, the differential of $\mu_K$ at $(z,w)$ is surjective if and only if $\{ \iota^* e^j ~|~ j \in J_{(z,w)} \}$ spans $\mathfrak{k}^*$. 
Since $\{ \iota^* e^0, \dots, \widehat{\iota^* e^k}, \dots, \iota^* e^n \}$ forms a basis  of $\mathfrak{k}^*$ for $k=0, \dots, n$, the lemma follows. \hfill $\Box$

If $(\alpha,\beta) \in \mathfrak{k}^* \oplus (\mathfrak{k}^* \!\otimes\! \mathbb{C})= \mathfrak{k}^* \otimes \mathrm{Im}\mathbb{H}$ is a regular value of $\mu_K$, then we have a smooth hyperK\"ahler quotient
\begin{align*}
M(\alpha,\beta)= (\mu_K^{-1}(\alpha,\beta)/K, g,I_1,I_2,I_3),
\end{align*}
on which $T^1=T^{n+1}/K$ acts, preserving the hyperK\"ahler structure. 
If $(\alpha,\beta)$ is a critical value, then the hyperK\"ahler quotient $M(\alpha,\beta)= \mu_K(\alpha,\beta)/K$ has singularities. 
A point in $M(\alpha,\beta)$, which is represented by $(z,w) \in \mu_K^{-1}(\alpha,\beta)$, is denoted by $[z,w]_K$. 

The hyperK\"ahler structure $(g,I_1,I_2,I_3)$ on $\mu_K^{-1}(\alpha,\beta)/K$, where $(\alpha,\beta)$ is a regular value of $\mu_K$, is defined as follows.
Consider the natural projection $$\pi \colon \mu_K^{-1}(\alpha,\beta) \to M(\alpha,\beta)$$ as a principal $K$-bundle. 
For $(z,w) \in \mu_K^{-1}(\alpha,\beta)$, the vertical subspace $V_{(z,w)}$ is defined to be the kernel of the differential 
\begin{align*}
\pi_{*(z,w)} \colon T_{(z,w)}\mu_K^{-1}(\alpha,\beta) \to T_{[z,w]_K}M(\alpha,\beta).
\end{align*}
Its orthogonal complement $H_{(z,w)}$ in $T_{(z,w)}\mu_K^{-1}(\alpha,\beta)$ is called the horizontal subspace.  
So we have the orthogonal decomposition 
\begin{align}\label{o-decomp}
T_{(z,w)}\mu_K^{-1}(\alpha,\beta)=H_{(z,w)} \oplus V_{(z,w)}.
\end{align}
It is easy to see that $V_{(z,w)}$ is the tangent space of the $K$-orbit through $(z,w)$ and that $H_{(z,w)}$ is a quaternionic subspace of $T_{(z,w)}\mathbb{H}^{n+1}$. 
The hyperK\"ahler structure on $M(\alpha,\beta)$ is defined so that the map $\pi_{*(z,w)}$ induces an hyperK\"ahler isometry from $H_{(z,w)}$ to $T_{[z,w]_K}M(\alpha,\beta)$. 

Define the action of $\mathbb{Z}_{n+1}= \{ \gamma \in \mathbb{C}~|~ \gamma^{n+1}=1 \}$ on $\mathbb{C}^2$ by
\begin{align*}
(u,v)\gamma=(u \gamma, v \gamma^{-1}).
\end{align*}
Since this action preserves the hyperK\"ahler structure of $\mathbb{C}^2$, 
the quotient space $\mathbb{C}^2/\mathbb{Z}_{n+1}$ is an orbifold which has the standard flat metric. 
A point in $\mathbb{C}^2/\mathbb{Z}_{n+1}$, which is represented by $(u,v) \in \mathbb{C}^2$, is denoted by $[u,v]_{\mathbb{Z}_{n+1}}$.
\begin{proposition}\label{prop-orb}
The map $\phi \colon \mathbb{C}^2/\mathbb{Z}_{n+1} \to M(0,0)$, which is defined by
\begin{align*}
\phi([u,v]_{\mathbb{Z}_{n+1}})
=[\frac{1}{\sqrt{n+1}}\begin{pmatrix}u \\ \vdots \\ u \end{pmatrix}, 
\frac{1}{\sqrt{n+1}}\begin{pmatrix}v \\ \vdots \\ v \end{pmatrix}]_K, 
\end{align*}
is an isometry. 
\end{proposition}
{\it Proof.}
Define the map $\widetilde{\phi} \colon \mathbb{C}^2 \to \mathbb{H}^{n+1}$ by 
\begin{align*}
\widetilde{\phi}(u,v)
=(\frac{1}{\sqrt{n+1}}\begin{pmatrix}u \\ \vdots \\ u \end{pmatrix}, 
\frac{1}{\sqrt{n+1}}\begin{pmatrix}v \\ \vdots \\ v \end{pmatrix}).
\end{align*}
Note that $(z,w) \in \mu_K^{-1}(0,0)$ if and only if 
\begin{align}\label{0-0}
|z_0|^2-|w_0|^2= \cdots = |z_n|^2-|w_n|^2, \quad z_0w_0= \cdots = z_n w_n.
\end{align}
Therefore, we see that $\widetilde{\phi}(u,v) \in \mu_K^{-1}(0,0)$ for each $(u,v) \in \mathbb{C}^2$.
Moreover, it is easy to see that $[\widetilde{\phi}(u,v)]_K=[\widetilde{\phi}(u^\prime,v^\prime)]_K$ if and only if there exists $\gamma \in \mathbb{Z}_{n+1}$ such that $(u^\prime, v^\prime)=(u,v)\gamma$. 
This implies that the map $\phi$ is well-defined and injective. 

To see that $\phi$ is surjective, fix an arbitrary $(z,w) \in \mu_K^{-1}(0,0)$.
By (\ref{0-0}), there exists $\zeta_k \in \mathbb{C}$ such that
\begin{align*}
|\zeta_k|=1, \quad (z_0,w_0)=(z_k \zeta_k, w_k \zeta_k^{-1}) \quad \text{for $k=1, \dots, n$.}
\end{align*}
If we set $\xi=(\zeta_1\cdots \zeta_n)^{\frac{-1}{n+1}}$, then we have $\zeta=(\xi, \zeta_1 \xi, \dots, \zeta_n \xi) \in K$ and 
\begin{align*}
[z,w]_K=[(z,w)\zeta]_K= \phi([z_0 \xi, w_0 \xi^{-1}]_{\mathbb{Z}_{n+1}}) .
\end{align*}
Thus we see that $\phi$ is surjective.

The differential $\widetilde{\phi}_{*(u,v)} \colon T_{(u,v)}\mathbb{C}^2 \to T_{\widetilde{\phi}(u,v)}\mathbb{H}^{n+1}$ is an isometric embedding for $(u,v) \in \mathbb{C}^2$. 
It is easy to see that the image of $\widetilde{\phi}_{*(u,v)}$ is perpendicular to the tangent space of the $K$-orbit through $\widetilde{\phi}(u,v)$ for $(u,v) \in \mathbb{C}^2 \setminus \{(0,0) \}$. 
This implies that $\widetilde{\phi}_{*(u,v)}$ induces an isometry from $T_{(u,v)}\mathbb{C}^2 $ to the horizontal subspace $H_{\widetilde{\phi}(u,v)}$ in (\ref{o-decomp}) for $(u,v) \in \mathbb{C}^2 \setminus \{(0,0) \}$. 
So we see that $\phi$ is an isometry. 
Thus we finish the proof. \hfill $\Box$

Thus we see that $(M(0,0),I_1)$ is isomorphic to $\mathbb{C}^2/\mathbb{Z}_{n+1}$
Moreover, it is well-known that $(M(\alpha,\beta), I_1)$ is a minimal resolution of $(M(0,0),I_1)$. 
\subsection{Torus actions on Ricci-flat ALE spaces of type $A_n$}
From now on, we fix
\begin{align}\label{alpha}
\alpha = \sum_{i=1}^n \alpha_i f^i \in \mathfrak{k}^*, \quad \text{where $\alpha_i>0$ for $i=1, \dots, n$}
\end{align}
and  $h \in (\mathfrak{t}^{n+1})^*$ such that $\alpha= \iota^* h$, that is,
\begin{align}\label{h}
h = \sum_{i=0}^n h_i e^i , \quad \text{where $h_0 \in \mathbb{R}$ arbitrary, $h_i=h_{i-1}+\alpha_i$ for $i=1, \dots, n$}.
\end{align}
If we set $h_{-1}=- \infty$ and $h_{n+1}=\infty$, then, by (\ref{alpha}) and (\ref{h}),  we have
\begin{align}\label{h2}
- \infty=h_{-1} <h_0 < h_1 < \dots < h_n < h_{n+1}=\infty.
\end{align}

By Proposition \ref{prop-rv}, $(\alpha,0) \in \mathfrak{k}^* \oplus (\mathfrak{k}^* \otimes \mathbb{C})= \mathfrak{k}^* \otimes \mathrm{Im}\mathbb{H}$ is a regular value of $\mu_K$.
Therefore $M(\alpha,0)= (\mu_K^{-1}(\alpha,0)/K, g,I_1,I_2,I_3)$ is a smooth manifold, on which $T^1=T^{n+1}/K$ acts, preserving the hyperK\"ahler structure. 
The K\"ahler form corresponding to $I_j$ is denoted by $\omega_j$ for $j=1,2,3$.
Then $\omega_\mathbb{C}=\omega_2+\sqrt{-1}\omega_3$ is a holomorphic symplectic form on $(M(\alpha,0),I_1)$.
Note that $(z,w) \in \mu_K^{-1}(\alpha,0)$ if and only if 
\begin{align}\label{alpha-0}
&\frac{1}{2}(|z_0|^2-|w_0|^2)-h_0= \dots =\frac{1}{2}(|z_n|^2-|w_n|^2)-h_n, \nonumber \\ 
&z_0 w_0= \dots =z_n w_n.
\end{align}

Define the action of a torus $G= \{ (\gamma_0,\gamma_1) \!\in\! \mathbb{C}^2~|~ |\gamma_0 |=|\gamma_1|=1 \}$ on $M(\alpha,0)$ by
\begin{align}\label{G-action}
[ \begin{pmatrix}z_0 \\z_1\\ \vdots \\z_n\end{pmatrix}, 
\begin{pmatrix}w_0 \\w_1\\ \vdots \\w_n\end{pmatrix}]_K (\gamma_0,\gamma_1)
=[ \begin{pmatrix}z_0\gamma_0\gamma_1 \\z_1\gamma_0\\ \vdots \\z_n\gamma_0\end{pmatrix}, 
\begin{pmatrix}w_0\gamma_1^{-1} \\w_1\\ \vdots \\w_n\end{pmatrix}]_K .
\end{align}
It is easy to see that the $G$-action is well-defined.
Note that the action of $(1, \gamma_1) \in G$ is the same as the action of $\rho(\gamma_1,1,\dots,1) \in T^1$, where $\rho \colon T^{n+1} \to T^1$ is the group homomorphism in (\ref{exact}). 
Denote the standard basis of $\mathfrak{g}$ by $p_0,p_1$ and the dual basis by $p^0,p^1 \in \mathfrak{g}^*$.
Then we have the following. 
Since the proof is straightforward, we omit the proof.
\begin{proposition}\label{prop-moment}
$(1)$ The $G$-action on $M(\alpha,0)$ preserves the complex structure $I_1$ and the corresponding K\"ahler form $\omega_1$.
Moreover, this action admits a moment map $\mu_G \colon M(\alpha,0) \to \mathfrak{g}^*$, which is given by
\begin{align}\label{moment}
\mu_G([z,w]_K)=\frac{1}{2}(\sum_{i=0}^n |z_i|^2)p^0 + \{ \frac{1}{2}(|z_k|^2-|w_k|^2)-h_k \} p^1.
\end{align}
$(2)$ $R_{\mathrm{Exp}_G\xi}^*\omega_\mathbb{C}=e^{\sqrt{-1}\langle p^0, \xi \rangle}\omega_\mathbb{C}$ holds for $\xi \in \mathfrak{g}$, where $\mathrm{Exp}_G \colon \mathfrak{g} \to G$ is the exponential map for $G$.
\end{proposition}
We remark that, in (\ref{moment}), the term $\frac{1}{2}(|z_k|^2-|w_k|^2)-h_k$ is independent of  $k=0,1, \dots ,n$ due to (\ref{alpha-0}). 

We note that, due to (\ref{h2}),  for each $y \in \mathbb{R}$, there exists unique $k_0 \in \{ 0, \dots, n+1 \}$ such that $-h_{k_0 -1} >y \ge -h_{k_0}$ holds. 
\begin{proposition}\label{prop-image}
$(1)$ The image $\mathrm{Im}\mu_G$ of the moment map $\mu_G \colon M(\alpha,0) \to \mathfrak{g}^*$ is given by 
\begin{align}\label{image}
\mathrm{Im}\mu_G 
= \{ xp^0 + yp^1 \in \mathfrak{g}^*~|~ x \ge \sum_{i=k_0}^n(y+ h_i)\quad\text{if  $-h_{k_0 -1} >y \ge -h_{k_0}$}\}.
\end{align}
$(2)$ For each $xp^0 + yp^1 \in \mathrm{Im}\mu_G$, $\mu_G^{-1}(xp^0 + yp^1)$ consists of a single $G$-orbit in $M(\alpha,0)$. 
\end{proposition}
{\it Proof.}
$(1)$ Denote the right hand side of (\ref{image}) by $\Delta$. 

First we show $\mathrm{Im}\mu_G \subset \Delta$.
Fix any $xp^0 + yp^1=\mu_G([z,w]_K) \in \mathrm{Im}\mu_G$. 
Then there exists unique $k_0 \in \{ 0, \dots, n+1 \}$ such that $-h_{k_0 -1} >y \ge -h_{k_0}$ holds. 
Due to (\ref{alpha-0}), there exists $d_0  \in \mathbb{C}$ such that  
\begin{align}\label{zw2}
\frac{1}{2}(|z_i|^2-|w_i|^2)-h_i=y, \quad z_i w_i=d_0 \quad \text{for $i=0, \dots, n$}.
\end{align}
Since $\{ \frac{1}{2}(|z_i|^2+|w_i|^2) \}^2= \{ \frac{1}{2}(|z_i|^2-|w_i|^2) \}^2 + |z_i w_i|^2=(y+ h_i)^2+|d_0|^2$, we have 
\begin{align}\label{zw}
|z_i|^2=\sqrt{(y+ h_i)^2+|d_0|^2}+ (y+ h_i), \quad 
|w_i|^2=\sqrt{(y+ h_i)^2+|d_0|^2}- (y+ h_i). 
\end{align}
So we have 
\begin{align}
x = \frac{1}{2}\sum_{i=0}^n |z_i|^2 &= \frac{1}{2}\sum_{i=0}^n\{ \sqrt{(y+ h_i)^2+|d_0|^2} +(y+ h_i) \} \label{ineq}\\
&\ge \frac{1}{2}\sum_{i=0}^n\{ |y+ h_i| +(y+ h_i) \} =\sum_{i=k_0}^n(y+ h_i).
\end{align}
Thus we see $\mathrm{Im}\mu_G \subset \Delta$.

Next we show $\mathrm{Im}\mu_G \supset \Delta$.
Fix any $xp^0 + yp^1 \in \Delta$. 
Then there exists unique $k_0 \in \{ 0, \dots, n+1 \}$ such that $-h_{k_0 -1} >y \ge -h_{k_0}$ holds. 
We want to show that there exists $[z,w]_K \in M(\alpha,0)$ such that $\mu_G([z,w]_K)=xp^0 +yp^1$. 
By (\ref{ineq}),  it is equivalent to show that there exist $(z,w) \in \mathbb{C}^{n+1} \times \mathbb{C}^{n+1}$ and $d_0  \in \mathbb{C}$ such that (\ref{zw2})  holds and $f_y(d_0)=x$, where $f_y \colon \mathbb{C} \to \mathbb{R}$ is a function defined by
\begin{align}\label{def-f} 
f_y(d)= \frac{1}{2}\sum_{i=0}^n\{ \sqrt{(y+ h_i)^2+|d|^2} +(y+ h_i) \}.
\end{align}

On the other hand, since $xp^0 + yp^1 \in \Delta$, we have 
\begin{align}\label{prop-f} 
f_y(0)= \frac{1}{2}\sum_{i=0}^n\{ |y+ h_i| +(y+ h_i) \} =\sum_{i=k_0}^n(y+ h_i) \le x.
\end{align}
So there exists $d_0 \in \mathbb{C}$ such that $f_y(d_0)=x$. 
Moreover, by (\ref{zw}), there exists $(z,w) \in \mathbb{C}^{n+1} \times \mathbb{C}^{n+1}$ such that (\ref{zw2})  holds.  
So we see that there exists $[z,w]_K \in M(\alpha,0)$ such that $\mu_G([z,w]_K)=xp^0 +yp^1$. 
Thus we see that $\mathrm{Im}\mu_G \supset \Delta$.
\\
$(2)$ Fix any $xp^0 + yp^1 \in \mathrm{Im}\mu_G$. 
In the proof of $(1)$, we constructed $[z,w]_K \in \mu_G^{-1}(xp^0 + yp^1)$. 
In the construction, there are some ambiguities for the choice of $d_0 \in \mathbb{C}$ and $(z,w) \in \mathbb{C}^{n+1} \times \mathbb{C}^{n+1}$.
It is easy to see that the ambiguity for $[z,w]_K$ corresponds to the $G$-action. 
\hfill$\Box$

Set 
\begin{align}\label{equ-lk}
l_k= \{ xp^0 + yp^1| -h_{k-1} \ge y \ge -h_k, ~x= \sum_{i=k}^n(y+ h_i) \}
\end{align}
 for $k=0,\dots, n+1$, and set $v_k=l_k \cap l_{k+1}$ for $k=0,\dots, n$.
Then the image of the moment map $\mu_G$ is decomposed in the following way (Fig.1):
\begin{align}\label{decomposition}
\mathrm{Im} \mu_G = \mathrm{int}(\mathrm{Im} \mu_G) \cup ~\bigcup_{k=0}^{n+1}\mathrm{int}(l_k) ~ \cup \{ v_0 ,\dots. v_n \},
\end{align}
where $\mathrm{int}(\mathrm{Im} \mu_G)$ and $\mathrm{int}(l_k)$ are the interior of $\mathrm{Im} \mu_G$ and $l_k$, respectively. 
Note that $\bigcup_{k=1}^n \mu_G^{-1}(l_k)$ is the exceptional divisor for the minimal resolution of $\mathbb{C}^2/\mathbb{Z}_{n+1}$.

Next we determine the isotropy subgroup $G_{[z,w]_K}$ of $G$ at $[z,w]_K \in M(\alpha,0)$ as follows. 
\begin{proposition}\label{prop-isotropy}
$(1)$ If $\mu_G([z,w]_K) \in \mathrm{int}(\mathrm{Im} \mu_G)$, then $G_{[z,w]_K}$ is the trivial subgroup and $z_i \ne 0$, $w_i \ne 0$ for $i=0,\dots, n$. 
\\
$(2)$ If $\mu_G([z,w]_K) \in \mathrm{int}(l_{k_0})$, then $G_{[z,w]_K}$ is the subtorus $H_{1,-(n+1-k_0)}$, whose Lie algebra is generated by $p_0 -(n+1-k_0)p_1$. 
Moreover, the following holds.
\begin{align}\label{zw3} 
&z_i=0 \quad\text{for $i=0, \dots, k_0-1$},\quad
w_i \ne 0 \quad\text{for $i=0, \dots, k_0-1$}\nonumber \\
&z_i \ne 0 \quad\text{for $i=k_0, \dots, n$},\hspace{9mm}
w_i=0 \quad\text{for $i=k_0, \dots, n$}.
\end{align}
$(3)$ If $\mu_G([z,w]_K) =v_{k_0}$, then $G_{[z,w]_K}=G$ and
\begin{align}\label{zw4}
&z_i=0 \quad\text{for $i=0, \dots, k_0$},\hspace{9.7mm}
w_i \ne 0 \quad\text{for $i=0, \dots, k_0-1$}\nonumber \\
&z_i \ne 0 \quad\text{for $i=k_0+1, \dots, n$},\quad
w_i=0 \quad\text{for $i=k_0, \dots, n$}. 
\end{align}
In particular, $\mu_G^{-1}(v_k)$ is a single point in $M(\alpha,0)$ for $k=0,\dots, n$.
\end{proposition}
{\it Proof.}
By the proof of Proposition \ref{prop-image}, the condition $xp^0+yp^1 \in \mathrm{Im} \mu_G$ holds if and only if there exist $(z,w) \in \mathbb{C}^{n+1} \times \mathbb{C}^{n+1}$ and $d_0  \in \mathbb{C}$ such that (\ref{zw2}) holds and $f_y(d_0)=x$, where the function $f_y$ is defined by (\ref{def-f}).  
\\
$(1)$ $\mu_G([z,w]_K) \in \mathrm{int}(\mathrm{Im} \mu_G)$ if and only if there exists $k_0 \in \{ 0, \dots, n+1\}$ such that $-h_{k_0 -1} >y \ge -h_{k_0}$ and $x > \sum_{i=k_0}^n(y+ h_i)$, where $xp^0+yp^1=\mu_G([z,w]_K)$. 
Since $x > \sum_{i=k_0}^n(y+ h_i)$, (\ref{prop-f}) implies $d_0 \ne 0$. 
Then (\ref{zw2}) implies the claim. 
\\
$(2)$ $\mu_G([z,w]_K) \in \mathrm{int}(l_{k_0})$ if and only if $-h_{k_0 -1} >y > -h_{k_0}$ and $x = \sum_{i=k_0}^n(y+ h_i)$, where $xp^0+yp^1=\mu_G([z,w]_K)$. 
Since $x = \sum_{i=k_0}^n(y+ h_i)$, (\ref{prop-f}) implies $d_0 = 0$. 
Then (\ref{zw2}) implies (\ref{zw3}). 
Moreover, the isotropy subgroup can be easily determined by the following. 
\begin{align*}
[\begin{pmatrix} 0 \\ \vdots \\ 0 \\ z_{k_0} \\ z_{k_0 +1} \\ \vdots \\ z_n \end{pmatrix}\!\!,\begin{pmatrix} w_0  \\ \vdots \\ w_{k_0 -1} \\ 0 \\ 0 \\ \vdots \\ 0\end{pmatrix}]_K
(e^{\sqrt{-1}s}, e^{\sqrt{-1}t})
=[\begin{pmatrix} 0 \\ \vdots \\ 0 \\ z_{k_0}e^{\sqrt{-1}\{(n+1-k_0)s+ t \}} \\ z_{k_0 +1} \\ \vdots \\ z_n \end{pmatrix}\!\!,\begin{pmatrix} w_0  \\ \vdots \\ w_{k_0 -1} \\ 0 \\ 0 \\ \vdots \\ 0\end{pmatrix}]_K. 
\end{align*}
$(3)$ $\mu_G([z,w]_K) =v_{k_0}$ if and only if $y = -h_{k_0}$ and $x = \sum_{i=k_0}^n(y+ h_i)$, where $xp^0+yp^1=\mu_G([z,w]_K)$. 
By the same argument as in the proof of $(2)$, we see the claim.  \hfill $\Box$

Next we define a good local coordinate around a fixed point. 
Let $v_{k_0} \in \mathrm{Im}\mu_G$ as in (\ref{decomposition}) for $k_0 =0, 1, \dots, n$. 
Since $\mu_G^{-1}(v_{k_0})$ is a single point, we set $P_{k_0}=\mu_G^{-1}(v_{k_0})$. 
Define an open set $U_{k_0}$ in $M(\alpha,0)$ by
\begin{align*}
U_{k_0}= \{ [z,w]_K ~|~ \text{$w_i \ne 0$ for $i=0, \dots, k_0 -1$, $z_j \ne 0$ for $j=k_0 +1, \dots, n$}\}.
\end{align*}
Define a map $\varphi_{k_0} \colon U_{k_0} \to \mathbb{C}^2$ by
\begin{align}\label{equ-coord}
\varphi_{k_0}([z,w]_K)= \Big( 
z_{k_0} \prod_{i=0}^{k_0 -1}&\frac{\sqrt{2(h_{k_0}-h_i)}}{w_i}\prod_{j=k_0 +1}^{n}\frac{z_j}{\sqrt{2(h_j-h_{k_0})}}, \nonumber \\
& w_{k_0} \prod_{i=0}^{k_0 -1}\frac{w_i}{\sqrt{2(h_{k_0}-h_i)}}\prod_{j=k_0 +1}^{n}\frac{\sqrt{2(h_j-h_{k_0})}}{z_j}\Big). 
\end{align}
It is easy to see that $\varphi_{k_0} \colon U_{k_0} \to \mathbb{C}^2$ is well-defined and that $P_{k_0} \in U_{k_0}$ and $\varphi_{k_0}(P_{k_0})=(0,0)$.
Moreover, we have the following.
\begin{proposition}\label{prop-coord}
$(1)$ $\varphi_{k_0} \colon U_{k_0} \to \mathbb{C}^2$ is a holomorphic local coordinate of $M(\alpha,0)$ with respect to $I_1$.
\\
$(2)$ If we write $\varphi_{k_0}(q)=(u_1(q),u_2(q))$ for $q \in U_{k_0}$ and $u_i=x_i +  \sqrt{-1}y_i$, where $x_i, y_i \in \mathbb{R}$, for $i=1,2$, then $(\frac{\partial}{\partial x_1})_{P_{k_0}}, (\frac{\partial}{\partial y_1})_{P_{k_0}}, (\frac{\partial}{\partial x_2})_{P_{k_0}}, (\frac{\partial}{\partial y_2})_{P_{k_0}}$ is an orthonormal basis of $T_{P_{k_0}}(M(\alpha,0))$.
\end{proposition}
{\it Proof.}
$(1)$ To prove the claim, we need to describe $(M(\alpha,0), I_1)$ as a quotient in geometric invariant theory. 
Note that $M(\alpha,0)$ is an example of toric hyperK\"ahler varieties, which were introduced in \cite{BD}. 
We refer the readers to \cite{Ko} for basic properties of toric hyperK\"ahler varieties. 

Define the action of $T^{n+1}_\mathbb{C} = \{ \zeta=(\zeta_0, \cdots, \zeta_n) ~|~ \zeta_i  \in \mathbb{C}^\times \quad \text{for $i=0, \dots, n$} \}$, where $\mathbb{C}^\times = \mathbb{C} \setminus \{ 0 \}$, on $\mathbb{H}^{n+1} = \mathbb{C}^{n+1} \times \mathbb{C}^{n+1}$ by (\ref{equ-action}).
Define a group homomorphism $\rho_\mathbb{C} \colon T^{n+1}_\mathbb{C} \to \mathbb{C}^\times$ by $\rho_\mathbb{C}(\zeta)=\zeta_0 \dots \zeta_n$.
Let $K_\mathbb{C}$ be the kernel of $\rho_\mathbb{C}$, which is the complexification of the torus $K$. 
If we set 
$$\mu_{K, \mathbb{C}}=\mu_K^2+\sqrt{-1}\mu_K^3 \colon \mathbb{C}^{n+1} \times \mathbb{C}^{n+1} \to \mathfrak{k}_\mathbb{C},~~\text{where $\mathfrak{k}_\mathbb{C}=\mathfrak{k} \otimes _\mathbb{R} \mathbb{C}$, }$$
then we have
\begin{align*}
\mu_{K, \mathbb{C}}^{-1}(0)= \{ (z,w) \in \mathbb{C}^{n+1} \times \mathbb{C}^{n+1}~|~ z_0w_0= \dots=z_n w_n \}.
\end{align*}
So $\mu_{K, \mathbb{C}}^{-1}(0)$ is an algebraic variety, on which $K_\mathbb{C}$ acts. 
Then we have a notion of stability for the action of $K_\mathbb{C}$ on $\mu_{K, \mathbb{C}}^{-1}(0)$ in geometric invariant theory according to $\alpha \in \mathfrak{k}$, which we have fixed in (\ref{alpha}). 
Since the torus $K$ acts on $\mu_K^{-1}(\alpha,0)$ freely in our case, 
$\alpha$-semistability is equivalent to $\alpha$-stability. 
So, if we denote the set of $\alpha$-stable points by $\mu_{K, \mathbb{C}}^{-1}(0)^{\alpha-st}$, then we have
\begin{align*}
(M(\alpha,0), I_1)=\mu_{K, \mathbb{C}}^{-1}(0)^{\alpha-st}/K_\mathbb{C},
\end{align*}
where the right hand side is a geometric quotient. 
A point in $\mu_{K, \mathbb{C}}^{-1}(0)^{\alpha-st}/K_\mathbb{C}$, which is represented by $(z,w) \in \mu_{K, \mathbb{C}}^{-1}(0)^{\alpha-st}$, is denoted by $[z,w]_{K_\mathbb{C}}$. 
\begin{claim}
Set
\begin{align*}
\tilde{U}_{k_0}\!=\! \{ (z,w) \in \mu_{K, \mathbb{C}}^{-1}(0) ~|~ \text{$w_i \!\ne\! 0$ for $i=0, \dots, k_0 \!-\!1$, $z_j \!\ne\! 0$ for $j=k_0 \!+\!1, \dots, n$}\}.
\end{align*}
Then $\tilde{U}_{k_0} \subset \mu_{K, \mathbb{C}}^{-1}(0)^{\alpha-st}$ holds.
\end{claim}
{\it Proof.}
Due to Lemma 3.6 in \cite{Ko}, a point $(z,w) \in \mu_{K, \mathbb{C}}^{-1}(0)$ is $\alpha$-stable if and only if
\begin{align}\label{equ-stable}
\alpha \in \sum_{j=0}^n \mathbb{R}_{>0}|z_j|^2 \iota^* e^j - \sum_{i=0}^n \mathbb{R}_{>0}|w_i|^2 \iota^* e^i.
\end{align}
On the other hand, by (\ref{alpha}) and  (\ref{equ-f}), we have 
\begin{align}\label{equ-alpha2}
\alpha \in \sum_{j=k_0 +1}^n \mathbb{R}_{>0} \iota^* e^j - \sum_{i=0}^{k_0 -1} \mathbb{R}_{>0} \iota^* e^i.
\end{align}
So the claim follows from (\ref{equ-stable}) and (\ref{equ-alpha2}). \hfill $\Box$

So we have $\tilde{U}_{k_0} / K_\mathbb{C} = U_{k_0}$. 
Moreover, $\varphi_{k_0}([z,w]_{K_\mathbb{C}})$ is also given by the right hand side of (\ref{equ-coord}) for $(z,w) \in \tilde{U}_{k_0}$. 
So $\varphi_{k_0} \colon U_{k_0} \to \mathbb{C}^2$ is holomorphic with respect to $I_1$.
If we write $\varphi_{k_0}([z,w]_{K_\mathbb{C}})=(u_1(z,w),u_2(z,w))$ for $(z,w) \in \tilde{U}_{k_0}$, then we have
\begin{align}\label{equ-holo}
[z,w]_{K_\mathbb{C}}
=
[\begin{pmatrix} \vdots \\ \frac{z_i w_i}{\sqrt{2(h_{k_0}-h_i)}} \\ \vdots \\ u_1(z,w) \\ \vdots \\ \sqrt{2(h_j-h_{k_0})} \\ \vdots \end{pmatrix},\begin{pmatrix} \vdots  \\ \sqrt{2(h_{k_0}-h_i)} \\ \vdots \\ u_2(z,w) \\ \vdots \\ \frac{z_j w_j}{\sqrt{2(h_j-h_{k_0})}} \\ \vdots \end{pmatrix}]_{K_\mathbb{C}},
\end{align}
where $u_1(z,w), u_2(z,w)$ in the right hand side are the $k_0$-th component. 
Since $z_0 w_0 = \dots = z_n w_n = u_1(z,w)u_2(z,w)$, the right hand side in (\ref{equ-holo}) depends only on $u_1(z,w),~u_2(z,w)$. 
So it can be seen that $ \varphi_{k_0} \colon U_{k_0} \to \mathbb{C}^2$ is bijective, moreover, a homeomorphism. 
Thus we see that $\varphi_{k_0} \colon U_{k_0} \to \mathbb{C}^2$ is a holomorphic local coordinate.
\\
$(2)$ Define smooth curves $\tilde{\gamma}_i \colon (-\epsilon, \epsilon) \to \mu_K^{-1}(\alpha,0)$ for $i=1,2$ by
\begin{align*}
&\tilde{\gamma}_1(t)
=
\Big(\begin{pmatrix} \vdots \\ 0 \\ \vdots \\ t \\ \vdots \\ \sqrt{2(h_j-h_{k_0})+t^2} \\ \vdots \end{pmatrix},\begin{pmatrix} \vdots  \\ \sqrt{2(h_{k_0}-h_i)-t^2} \\ \vdots \\ 0 \\ \vdots \\ 0 \\ \vdots \end{pmatrix}\Big), \\
&\tilde{\gamma}_2(t)
=
\Big(\begin{pmatrix} \vdots \\ 0 \\ \vdots \\ 0 \\ \vdots \\ \sqrt{2(h_j-h_{k_0})-t^2} \\ \vdots \end{pmatrix},\begin{pmatrix} \vdots  \\ \sqrt{2(h_{k_0}-h_i)+t^2} \\ \vdots \\ t \\ \vdots \\ 0 \\ \vdots \end{pmatrix}\Big),
\end{align*}
for small $\epsilon > 0$.
Then $\tilde{\gamma}_i$ induces smooth curves $\gamma_i \colon (-\epsilon, \epsilon) \to M(\alpha,0)$ for $i=1,2$. 
Then we have $\gamma_1(0)=\gamma_2(0)=P_{k_0}$ and $\frac{d(\varphi_{k_0} \circ \gamma_i)}{dt}(0)=(\frac{\partial}{\partial x_i})_{P_{k_0}}$.
Moreover, $\frac{d \tilde{\gamma}_1}{dt}(0)$ is perpendicular to $\frac{d \tilde{\gamma}_2}{dt}(0)$, and $\frac{d \tilde{\gamma}_i}{dt}(0)$ are also perpendicular to the $K$-orbit through $\tilde{\gamma}_1(0)=\tilde{\gamma}_2(0)$. 
Then the claim follows. \hfill $\Box$
\subsection{Involutions on Ricci-flat ALE spaces of type $A_n$}\label{sec-involution}
Define a map $\sigma \colon M(\alpha,0) \to M(\alpha,0)$ by 
\begin{align*}
\sigma([z,w]_K)=[\overline{z},\overline{w}]_K.
\end{align*}
It is easy to see that this map is well-defined. 
Then the following holds.
\begin{proposition}\label{prop-involution}
$(1)$ $\sigma \colon M(\alpha,0) \to M(\alpha,0)$ is an isometric involution, which is anti-holomorphic with respect to $I_1$ and holomorphic with respect to $I_2$. 
\\
$(2)$ $\sigma([z,w]_K (\gamma_0,\gamma_1))=\sigma([z,w]_K) (\gamma_0,\gamma_1)^{-1}$ holds for any $[z,w]_K \in M(\alpha,0)$ and $(\gamma_0,\gamma_1) \in G$.
\end{proposition}
{\it Proof.}
$(1)$ Define a map $\widetilde{\sigma} \colon \mathbb{H}^{n+1} \to \mathbb{H}^{n+1}$ by $\widetilde{\sigma}(z,w)=(\overline{z},\overline{w})$ for $(z,w) \in \mathbb{C}^{n+1}\times \mathbb{C}^{n+1}=\mathbb{H}^{n+1}$. 
It is easy to see that $\widetilde{\sigma}$ is an isometric involution, which is anti-holomorphic with respect to $I_1$ and holomorphic with respect to $I_2$. 
Moreover, it is easy to see that the differential $\widetilde{\sigma}_{*(z,w)} \colon T_{(z,w)}\mathbb{H}^{n+1} \to T_{\widetilde{\sigma}(z,w)}\mathbb{H}^{n+1}$ induces  an isometry from the horizontal subspace $H_{(z,w)}$ to $H_{\widetilde{\sigma}(z,w)}$ in (\ref{o-decomp}) for $(z,w) \in \mu_K^{-1}(\alpha,0)$.
Therefore the induced map $\sigma$ is also an isometric involution, which is anti-holomorphic with respect to $I_1$ and holomorphic with respect to $I_2$.  
\\
$(2)$ is obvious. \hfill $\Box$

To describe the set of fixed points of the involution $\sigma \colon M(\alpha,0) \to M(\alpha,0)$, we set $G_\mathbb{R}=G \cap \mathbb{R}^2$, which consists of four elements $(1,1)$, $(-1,1)$, $(1,-1)$, $(-1,-1)$. 
\begin{proposition}\label{prop-fixedpts}
Let $M^\sigma$ be the set of fixed points of $\sigma$. 
Then the following holds.
\\
$(1)$ $M^\sigma= \{ \mu_K^{-1}(\alpha,0) \cap (\mathbb{R}^{n+1} \times \mathbb{R}^{n+1}) \}/(K \cap \mathbb{R}^{n+1}).$
\\
$(2)$ For any $xp^0 +yp^1 \in \mathrm{Im} \mu_G$, $\mu_G^{-1}(xp^0 +yp^1) \cap M^\sigma$ is a single $G_\mathbb{R}$-orbit. 
\\
$(3)$ Take a subset $M^\sigma_{++}$ of $M^\sigma$ so that $\mu_G \colon M \to \mathfrak{g}^*$ induces a homeomorphism from $M^\sigma_{++}$ to $\mathrm{Im} \mu_G$. 
Set 
\begin{align*}
&M^\sigma_{-+}= \{ [z,w]_K(-1, 1) ~|~ [z,w]_K \in M^\sigma_{++} \}, \\
&M^\sigma_{+-}= \{ [z,w]_K(1, -1) ~|~ [z,w]_K \in M^\sigma_{++} \}, \\
&M^\sigma_{--}= \{ [z,w]_K(-1, -1) ~|~ [z,w]_K \in M^\sigma_{++} \}.
\end{align*}
Then $M^\sigma = M^\sigma_{++} \cup M^\sigma_{-+} \cup M^\sigma_{+-} \cup M^\sigma_{--}$ holds. 
Moreover, for $k=0, \dots, n+1$, 
\begin{align*}
&\text{$M^\sigma_{++}$and $M^\sigma_{-+}$ are glued along $l_k$ if and only if $n-k$ is odd,}\\
&\text{$M^\sigma_{+-}$and $M^\sigma_{--}$ are glued along $l_k$ if and only if $n-k$ is odd,}\\
&\text{$M^\sigma_{++}$and $M^\sigma_{--}$ are glued along $l_k$ if and only if $n-k$ is even,}\\
&\text{$M^\sigma_{+-}$and $M^\sigma_{-+}$ are glued along $l_k$ if and only if $n-k$ is even.}
\end{align*}
Therefore, $M^\sigma$ is homeomorphic to 
\begin{align*}
&\text{a two-holed orientable surface of genus $\frac{n-1}{2}$ if $n$ is odd,}\\
&\text{a one-holed orientable surface of genus $\frac{n}{2}$ if $n$ is even.}
\end{align*}\end{proposition}
{\it Proof.}
$(1)$ Set $M_\mathbb{R}=\{ \mu_K^{-1}(\alpha,0) \cap (\mathbb{R}^{n+1} \times \mathbb{R}^{n+1}) \}/(K \cap \mathbb{R}^{n+1})$.
Since $M^\sigma \supset M_\mathbb{R}$ is obvious, we show that $M^\sigma \subset M_\mathbb{R}$.
Fix $[z,w]_K \in M^\sigma$. 
Then there exists $\zeta \in K$ such that $(\overline{z},\overline{w})=(z\zeta, w\zeta^{-1})$. 
So we have $(\overline{z}\zeta^{-\frac{1}{2}},\overline{w}\zeta^\frac{1}{2})=(z\zeta^\frac{1}{2},w\zeta^{-\frac{1}{2}})$. 
This implies that $(z\zeta^\frac{1}{2},w\zeta^{-\frac{1}{2}}) \in \mu_K^{-1}(\alpha,0) \cap (\mathbb{R}^{n+1} \times \mathbb{R}^{n+1})$.
Thus we see that $[z,w]_K = [z\zeta^\frac{1}{2},w\zeta^{-\frac{1}{2}}]_K \in M_\mathbb{R}$.
\\
$(2)$ Fix $xp^0 + yp^1 \in \mathrm{Im} \mu_G$. 
Since $\mu_G^{-1}(xp^0 + yp^1)$ is a single $G$-orbit by Proposition \ref{prop-image} $(2)$, it is enough to show that $\mu_G^{-1}(xp^0 + yp^1) \cap M^\sigma \ne \emptyset$.
As in the proof of Proposition \ref{prop-image} $(1)$, it is equivalent to show that there exist $(z,w) \in \mathbb{R}^{n+1} \times \mathbb{R}^{n+1}$ and $d_0  \in \mathbb{R}$ such that (\ref{zw2})  holds and $f_y(d_0)=x$, where $f_y \colon \mathbb{C} \to \mathbb{R}$ is a function defined by (\ref{def-f}). 
In fact, such $(z,w)$ and $d_0$ exist by the same argument as in the proof of Proposition \ref{prop-image} $(1)$. 
\\
$(3)$ By Proposition \ref{prop-isotropy}, for $k=0,\dots,n+1$, we have $G_{[z,w]_K}=H_{1,-(n+1-k)}$ if $\mu_G([z,w]_K) \in \mathrm{int}(l_k)$, where $\mathrm{int}(l_k)$ is the interior of $l_k$. 

On the other hand, for $k=0,\dots,n+1$, we have
\begin{align*}
& \text{$(-1,1) \in H_{1,-(n+1-k)}$ if and only if $n-k$ is odd,} \\
& \text{$(-1,-1) \in H_{1,-(n+1-k)}$ if and only if $n-k$ is even,} \\
& \text{$(1,-1) \not\in H_{1,-(n+1-k)}$ for any $k$.} 
\end{align*}
So the claim follows. \hfill $\Box$
\subsection{Lagrangian mean curvature flows in Ricci-flat ALE spaces of type $A_n$}
Suppose that $a,b \in \mathbb{Z}$ are coprime. 
Let $H_{a,b}$ be the subtorus of $G$, whose Lie algebra $\mathfrak{h}_{a,b}$ is generated by $ap_0+bp_1 \in \mathfrak{g}$.
Set $w_{a,b} =ap_0+bp_1 \in \mathfrak{h}_{a,b}$. 
Then there exists $w^{a,b} \in \mathfrak{h}_{a,b}^*$ such that $\langle w^{a,b}, w_{a,b} \rangle =1$.
Let $\iota_{a,b} \colon H_{a,b} \to G$ be the embedding, and $\iota_{a,b*} \colon \mathfrak{h}_{a,b} \to \mathfrak{g}$, $\iota_{a,b}^* \colon \mathfrak{g}^* \to \mathfrak{h}_{a,b}^*$ the induced maps. 
By Propositions \ref{prop-moment} and \ref{prop-involution}, we have the following. 
\begin{lemma}\label{lemma-hab}
$(1)$ $\iota_{a,b}^*p^0=aw^{a,b}$, $\iota_{a,b}^*p^1=b w^{a,b}$.
\\
$(2)$ The moment map $\mu_{H_{a,b}} = \iota_{a,b}^* \circ \mu_G\colon M(\alpha,0) \to \mathfrak{h}_{a,b}^*$ for the action of $H_{a,b}$ on $(M(\alpha,0), \omega_1)$ is given by
\begin{align*}
\mu_{H_{a,b}}([z,w]_K)=\{\frac{a}{2}(\sum_{i=0}^n |z_i|^2) + \frac{b}{2}(|z_k|^2-|w_k|^2)-b h_k \} w^{a,b}.
\end{align*}
\\
$(3)$ $R_{\mathrm{Exp}_{H_{a,b}}\xi}^*\omega_\mathbb{C}=e^{\sqrt{-1}a \langle w^{a,b}, \xi \rangle}\omega_\mathbb{C}$ holds for $\xi \in \mathfrak{\mathfrak{h}}_{a,b}$, where $\mathrm{Exp}_{H_{a,b}} \colon \mathfrak{h}_{a,b} \to H_{a,b}$ is the exponential map for $H_{a,b}$.
\\
$(4)$ $\sigma([z,w]_K h)=\sigma([z,w]_K) h^{-1}$ holds for any $[z,w]_K \in M(\alpha,0)$ and $h \in H_{a,b}$.
\end{lemma}
Denote the sets of fixed points of the action of $H_{a,b}$, $G$ on $M(\alpha,0)$ by $M(\alpha,0)^{H_{a,b}}$, $M(\alpha,0)^G$, respectively. 
Then we have the following.
\begin{lemma}\label{chi3}
Suppose that $a,b \in \mathbb{Z}$ are coprime and that $b \ne- la$ for $l=0,-1, \dots, \\-(n+1)$. 
Then the following holds.
\\
$(1)$ $M(\alpha,0)^{H_{a,b}}=M(\alpha,0)^G$ holds. 
In particular, $(w_{a,b})^\#_p \ne 0$ holds for each $p \in L=M^\sigma \setminus M(\alpha,0)^G$.
\\
$(2)$ Let $\chi$ be the vector field on $L=M^\sigma \setminus M(\alpha,0)^G$, which is defined in Theorem \ref{construction} $(4)$.
Then there exist constants $C_1 >0$ and $C_2 >0$ such that, if $p=[z,w]_K \in L$ satisfies $|p| \ge C_1$, then $|\chi_p| \le \frac{C_2}{|p|}$ holds, where $|p|= \sqrt{|z|^2 + |w|^2}$.
\end{lemma}
{\it Proof.} $(1)$ follows from Proposition \ref{prop-isotropy}.
\\
$(2)$ The vector field $\chi$ is defined by $\chi_p=I_p \tilde{\xi}(p)^\#_p$, where $\tilde{\xi}(p) \in \mathfrak{h}_{a,b}$ satisfies $g_p(\tilde{\xi}(p)^\#_p, \eta^\#_p)= \langle aw^{a,b}, \eta \rangle$ for any $\eta \in \mathfrak{h}_{a,b}$.
If we set $\tilde{\xi}(p)=\alpha_p w_{a,b}$, where $\alpha_p \in \mathbb{R}$, then we have $g_p(\alpha_p (w_{a,b})^\#_p, (w_{a,b})^\#_p)= a$. 
So we have $\alpha_p= \frac{a}{g_p((w_{a,b})^\#_p, (w_{a,b})^\#_p)}$.  
Thus we have $|\chi_p|= |I_p \tilde{\xi}(p)^\#_p|= |\alpha_p (w_{a,b})^\#_p|=\frac{|a|}{|(w_{a,b})^\#_p|}$. 
Therefore it is enough to show the following claim.
\begin{claim}\label{wab}
There exist constants $C_1 >0$ and $C_3 >0$ such that, if $p \in M(\alpha,0)$ satisfies $|p| \ge C_1$, then $|(w_{a,b})^\#_p| \ge C_3|p|$ holds.
\end{claim}
{\it Proof.}
Set $\widetilde{S}= \{ (z,w) \in \mathbb{C}^{n+1} \times \mathbb{C}^{n+1}~|~ |z|^2 + |w|^2=1 \}$. 
Since the action of $K$ on $ \mathbb{C}^{n+1} \times \mathbb{C}^{n+1}$ preserves $\widetilde{S}$, we can define $S(M(\alpha^\prime,0))=(\mu_K^{-1}(\alpha^\prime,0) \cap \widetilde{S})/K$, which is a compact subset of $M(\alpha^\prime,0)$ for each $\alpha^\prime \in \mathfrak{k}^*$. 
Since $b \ne 0, -(n+1)a$, it is easy to see that $(w_{a,b})^\#_p \ne 0$ for each $p \in S(M(0,0))$.  
So there exist a constant $C_3 >0$ and an open neighborhood $U_O$ of the origin $O \in \mathfrak{k}^*$ such that $|(w_{a,b})^\#_p| \ge C_3$ holds for each $\alpha^\prime \in U_O$ and $p \in S(M(\alpha^\prime,0))$. 

On the other hand, recall we have fixed $\alpha \in \mathfrak{k}$. 
Note that, if $p =[z,w]_K \in M(\alpha,0)$, then $\frac{p}{|p|}=[\frac{z}{|p|},\frac{w}{|p|}]_K \in S(M(\frac{\alpha}{|p|^2},0))$ and $|(w_{a,b})^\#_{\frac{p}{|p|}}|=\frac{1}{|p|}|(w_{a,b})^\#_p|$. 
There exists a constant $C_1 >0$ such that, if $p \in M(\alpha,0)$ satisfies $|p| \ge C_1$, then $\frac{\alpha}{|p|^2} \in U_O$. 
Then we have $\frac{1}{|p|}|(w_{a,b})^\#_p|=|(w_{a,b})^\#_{\frac{p}{|p|}}| \ge C_3$.
So Claim \ref{wab} follows. \hfill $\Box$
\\
Thus we finish the proof of Lemma \ref{chi3}. \hfill $\Box$

So we have the following theorem. 
\begin{theorem}\label{LMCFALE}
Suppose that $a,b \in \mathbb{Z}$ are coprime and that $b \ne- la$ for $l=0,-1, \dots, \\ -(n+1)$. 
Fix $c_0 \in \mathfrak{h}_{a,b}^*$ and set $c_t= c_0 - t a w^{a,b} \in \mathfrak{h}_{a,b}^*$. 
Suppose that $\mu_{H_{a,b}}^{-1}(c_t) \cap M(\alpha,0)^G = \emptyset$ for each $t \in [0,T)$.
Then the following holds.
\\
$(1)$ Set $V_{c_t} = \mu_{H_{a,b}}^{-1}(c_t) \cap L$, where $L=M^\sigma \setminus M(\alpha,0)^G$. 
Let $\chi$ be the vector field on $L$, which is defined in Theorem \ref{construction} $(4)$. 
Then, for each $p \in V_{c_0}$, the integral curve $\{ \gamma_p(t) \}$ of the vector field $\chi$ with $\gamma_p(0)=p$ exists for $t \in [0,T)$. 
Moreover, $\gamma_p(t) \in V_{c_t}$ holds.
\\
$(2)$ $\{ F_t \colon V_{c_0} \times H_{a,b} \to M(\alpha,0) \}_{t \in [0,T)}$, defined by $F_t(p,h)=\gamma_p(t)h$, is a Lagrangian mean curvature flow.  
\end{theorem}
{\it Proof.}
$(1)$ By the same arguments as in the proof of Proposition \ref{prop-existence}, the first claim follows from Lemma \ref{chi3}. 
The second claim follows from Theorem \ref{construction} $(5)$. 
\\
$(2)$ Theorem \ref{construction} $(5)$, together with Lemma \ref{lemma-hab} and Propositions \ref{construction2}, implies the claim. 
\hfill $\Box$

Let us describe the set $V_{c_0}=\mu_{H_{a,b}}^{-1}(c_0) \cap L$. 
If we set $r_{c_0}=(\iota_{a,b}^*)^{-1}(c_0) \cap \mathrm{Im}\mu_G$, then $r_{c_0}$ is a segment or a ray in $\mathfrak{g}^*$. 
$V_{c_0}$ is compact if and only if $r_{c_0}$ is bounded. 
Since $\mu_{H_{a,b}}^{-1}(c_0)= \mu_G^{-1}(r_{c_0})$,  we have
\begin{align}\label{equ-decomp2}
V_{c_0}&=\mu_{H_{a,b}}^{-1}(c_0) \cap L \nonumber \\ 
&= \mu_G^{-1}(r_{c_0}) \cap M^\sigma \nonumber \\
&= (\mu_G^{-1}(r_{c_0}) \cap M^\sigma_{++}) 
\cup (\mu_G^{-1}(r_{c_0}) \cap M^\sigma_{-+}) \nonumber \\
& \hspace*{20mm}\cup 
(\mu_G^{-1}(r_{c_0}) \cap M^\sigma_{+-}) \cup 
(\mu_G^{-1}(r_{c_0}) \cap M^\sigma_{--}) .
\end{align}
Thus $V_{c_0}$ is decomposed into four pieces 
\begin{align*}
\text{$\mu_G^{-1}(r_{c_0}) \cap M^\sigma_{++}$, $\mu_G^{-1}(r_{c_0}) \cap M^\sigma_{-+}$, $\mu_G^{-1}(r_{c_0}) \cap M^\sigma_{+-}$, $\mu_G^{-1}(r_{c_0}) \cap M^\sigma_{--}$.}
\end{align*}
Note that each piece is homeomorphic to the set $r_{c_0}$.
By Proposition \ref{prop-fixedpts} $(3)$, we see how these pieces are connected to each other. 
Thus we see whether $V_{c_0}$ is connected or not. 
Note that the four pieces in (\ref{equ-decomp2}) are transformed into themselves or other pieces by the action of $$G_\mathbb{R}= \{ (1,1), (1,-1),(-1,1), (-1,-1) \}. $$
Since we have assumed that $a,b \in \mathbb{Z}$ are coprime, it is easy to see that 
\begin{align}\label{hab}
& \text{$(1,-1) \in H_{a,b}$ if and only if $a$ is even and $b$ is odd,} \nonumber \\
& \text{$(-1,1) \in H_{a,b}$ if and only if $a$ is odd and $b$ is even,} \\
& \text{$(-1,-1) \in H_{a,b}$ if and only if $a$ is odd and $b$ is odd.} \nonumber 
\end{align}
So, in any case, if we set $V_{c_0}^0=V_{c_0} \cap \mu_G^{-1}(\mathrm{int}(\mathrm{Im}\mu_G))$, we see that the restriction $F_0|_{V_{c_0}^0 \times H_{a,b}} \colon V_{c_0}^0 \times H_{a,b} \to M(\alpha,0)$ is a two-to-one map onto its image. 
Thus we see that the map $F_t \colon V_{c_0} \times H_{a,b} \to M(\alpha,0)$ is generically a two-to-one map onto its image for any $t \in [0,T)$. 
 
Let us describe the above Lagrangian mean curvature flow $\{ F_t \colon V_{c_0} \times H_{a,b} \to M(\alpha,0) \}_{t \in [0,T)}$ in more detail. 
If $a=0$, the flow is static, that is, $F_t=F_0 \colon V_{c_0} \times H_{a,b} \to M(\alpha,0)$ for any $t \in [0,T)$. 
So we may assume $a>0$. 
We describe the flow in the following three cases:
\\
$(1)$ the case $a>0,~\frac{b}{a}>0$,
\\
$(2)$ the case $a>0,~0 > \frac{b}{a} > -(n+1),~ \frac{b}{a} \not\in \mathbb{Z}$, 
\\
$(3)$ the case $a>0,~-(n+1) > \frac{b}{a}$. 
\\
Since the third case is similar to the first one, we will discuss the first and second cases. 
\subsection{The case $a>0,~\frac{b}{a}>0$}\label{case1}
Recall the decomposition of $\mathrm{Im}\mu_G$ in (\ref{decomposition}). 
Define $t_i \in \mathbb{R}$ by $\iota_{a,b}^*(v_i)=c_0 - t_i a w^{a,b}(=c_{t_i})$ for $i=0, \dots, n$.  In the case $a>0,~\frac{b}{a}>0$, we have
$$
-\infty = t_{-1} < t_0< t_1 < \dots < t_n < t_{n+1}=\infty.
$$
Since $\mu_{H_{a,b}}^{-1}(c_0) \cap M(\alpha,0)^G = \emptyset$, there exists unique $k_0 \in \{ 0, 1, \dots, n+1 \}$ such that $t_{k_0 -1} < 0 < t_{k_0}$. 
Then $(\iota_{a,b}^*)^{-1}(c_0)$ intersects with $\mathrm{int}(l_{k_0})$. 
Moreover, the set $r_{c_0}=(\iota_{a,b}^*)^{-1}(c_0) \cap \mathrm{Im}\mu_G$ is a ray in $\mathfrak{g}^*$. 
Each piece $\mu_G^{-1}(r_{c_0}) \cap M^\sigma_{++}$, $\mu_G^{-1}(r_{c_0}) \cap M^\sigma_{-+}$, $\mu_G^{-1}(r_{c_0}) \cap M^\sigma_{+-}$ and $\mu_G^{-1}(r_{c_0}) \cap M^\sigma_{--}$ in (\ref{equ-decomp2}) is homeomorphic to the ray $r_{c_0}$.
By Proposition \ref{prop-fixedpts}, we see that $V_{c_0}=\mu_{H_{a,b}}^{-1}(c_0) \cap L$ consists of two connected components $V_{c_0}^{(1)}$,  $V_{c_0}^{(2)}$ as follows: 
\\
\hspace*{10mm}$(a)$ If $n-k_0$ is odd, 
\\
\hspace*{20mm}$V_{c_0}^{(1)}=(\mu_G^{-1}(r_{c_0}) \cap M^\sigma_{++}) 
\cup (\mu_G^{-1}(r_{c_0}) \cap M^\sigma_{-+})$, 
\\
\hspace*{20mm}$V_{c_0}^{(2)}=(\mu_G^{-1}(r_{c_0}) \cap M^\sigma_{+-}) \cup 
(\mu_G^{-1}(r_{c_0}) \cap M^\sigma_{--})$. 
\\
\hspace*{20mm}Each component is diffeomorphic to $\mathbb{R}$.
\\
\hspace*{10mm}$(b)$ If $n-k_0$ is even, 
\\
\hspace*{20mm}$V_{c_0}^{(1)}=(\mu_G^{-1}(r_{c_0}) \cap M^\sigma_{++}) 
\cup (\mu_G^{-1}(r_{c_0}) \cap M^\sigma_{--})$, 
\\
\hspace*{20mm}$V_{c_0}^{(2)}=(\mu_G^{-1}(r_{c_0}) \cap M^\sigma_{+-}) 
\cup (\mu_G^{-1}(r_{c_0}) \cap M^\sigma_{-+})$. 
\\
\hspace*{20mm}Each component is diffeomorphic to $\mathbb{R}$.

Note that $\mu_{H_{a,b}}^{-1}(c_t) \cap M(\alpha,0)^G = \emptyset$ for each $t \in [0,t_{k_0})$. 
By Theorem \ref{LMCFALE}, we have a Lagrangian mean curvature flow $\{ F_t \colon V_{c_0} \times H_{a,b} \to M(\alpha,0) \}_{t \in [0,t_{k_0})}$.  
By Proposition \ref{prop-fixedpts}, for $t \in (t_{k_0}, t_{k_0 +1} )$, the four pieces $\mu_G^{-1}(r_{c_t}) \cap M^\sigma_{++}$, $\mu_G^{-1}(r_{c_t}) \cap M^\sigma_{-+}$, $\mu_G^{-1}(r_{c_t}) \cap M^\sigma_{+-}$ and $\mu_G^{-1}(r_{c_t}) \cap M^\sigma_{--}$ are connected to each other in a different way from the case $t \in [0, t_{k_0})$. 
Therefore, it is impossible to extend the flow $\{ F_t \colon V_{c_0} \times H_{a,b} \to M(\alpha,0) \}_{t \in [0,t_{k_0})}$ continuously to $t \in [0,t_{k_0}+\epsilon)$ for any $\epsilon >0$. 
This implies that the Lagrangian mean curvature flow $\{ F_t \colon V_{c_0} \times H_{a,b} \to M(\alpha,0) \}_{t \in [0,t_{k_0})}$ develops a singularity at $P_{k_0}=\mu_G^{-1}(v_{k_0})$ when $t$ goes to $t_{k_0}$. 

Let us describe the singularity of the flow $\{ F_t  \}_{t \in [0,t_{k_0})}$. 
In Proposition \ref{prop-coord} we have constructed the holomorphic local coordinate $\varphi_{k_0} \colon U_{k_0} \to \mathbb{C}^2$ around $P_{k_0}$. 
We write $\varphi_{k_0}(q)=(u_1(q),u_2(q))$ for $q \in U_{k_0}$ and $u_i=x_i + \sqrt{-1}y_i$, where $x_i, y_i \in \mathbb{R}$, for $i=1,2$. 
By (\ref{equ-coord}), we have 
\begin{align*}
(u_1,u_2)\mathrm{Exp}_{H_{a,b}}s w_{a,b} 
& = (u_1 e^{\sqrt{-1} \lambda^{(k_0)}_1 s}, u_2 e^{\sqrt{-1} \lambda^{(k_0)}_2 s}), 
\end{align*}
where 
\begin{align}\label{lambda}
\lambda^{(k_0)}_1= a(n +1 -k_0)+b,~\lambda^{(k_0)}_2= -a(n-k_0)-b. 
\end{align}
Since we have assumed that $a>0$ and $\frac{b}{a}>0$, we have $\lambda^{(k_0)}_1 >0$ and $\lambda^{(k_0)}_2 <0$ for $k_0=0, \dots, n$.
Thus we have
\begin{align*}
(w_{a,b})^\# & = (\sqrt{-1} \lambda^{(k_0)}_1 u_1, \sqrt{-1} \lambda^{(k_0)}_2 u_2) \\
& = -\lambda^{(k_0)}_1 y_1 \frac{\partial}{\partial x_1} + \lambda^{(k_0)}_1 x_1 \frac{\partial}{\partial y_1}  -\lambda^{(k_0)}_2 y_2 \frac{\partial}{\partial x_2} + \lambda^{(k_0)}_2 x_2 \frac{\partial}{\partial y_2}.
\end{align*}
If we set $\mu=\langle \mu_{H_{a,b}}(\cdot), w_{a,b} \rangle \colon M(\alpha,0) \to \mathbb{R}$, then we have $i((w_{a,b})^\#)\omega_1=-d \mu$, where $\omega_1$ is the K\"ahler form of $(M(\alpha,0),I_1)$.  

Define $\omega_1^0, \omega_1^1 \in \Omega^2(U_{k_0})$ and $\mu^0,\mu^1 \colon U_{k_0} \to \mathbb{R}$ by
\begin{align*}
&\omega_1|_{U_{k_0}}=\omega_1^0 + \omega_1^1, ~~~~\omega_1^0 =dx_1 \wedge dy_1 + dx_2 \wedge dy_2, \\
&\mu|_{U_{k_0}}=\mu^0 + \mu^1, ~~~~\mu^0(u_1,u_2)=\mu(P_{k_0})+\frac{\lambda^{(k_0)}_1}{2}|u_1|^2 + \frac{\lambda^{(k_0)}_2}{2}|u_2|^2.
\end{align*}
By definition of $\mu^1$, we have $\mu^1(0,0)=0$. 
Since $i((w_{a,b})^\#)\omega_1^0=-d \mu^0$ holds, we have $i((w_{a,b})^\#)\omega_1^1=-d \mu^1$. 
Since $(w_{a,b})^\#_{(0,0)}=0$, we have $(d \mu^1)_{(0,0)}=0$. 
Moreover, by Proposition \ref{prop-coord} $(2)$, we have $(\omega_1^1)_{(0,0)}=0 $.
So we have 
\begin{align*}
\lim_{(u_1,u_2) \to (0,0)}\frac{|(d\mu^1)_{(u_1,u_2)}|}{\sqrt{|u_1|^2 + |u_2|^2}}
=\lim_{(u_1,u_2) \to (0,0)}\frac{|\{ i((w_{a,b})^\#)\omega_1^1 \}_{(u_1,u_2)}|}{\sqrt{|u_1|^2 + |u_2|^2}}
=0.
\end{align*}
Therefore we have
\begin{align}\label{equ-limmu}
\lim_{(u_1,u_2) \to (0,0)}\frac{\mu^1(u_1,u_2)}{|u_1|^2 + |u_2|^2}=0.
\end{align}

Since $c_t=c_0-taw^{a,b}=c_{k_0}+(t_{k_0}-t)aw^{a,b}$ holds, we have $\langle c_t, w_{a,b} \rangle = \mu(P_{k_0})+a(t_{k_0} -t)$. 
Then we have
\begin{align}\label{equ-local}
(u_1,u_2) &\in \mu_{H_{a,b}}^{-1}(c_t) \cap U_{k_0} \nonumber \\ 
&\Longleftrightarrow \mu(u_1,u_2)=\mu(P_{k_0})+a(t_{k_0}-t) \nonumber \\
&\Longleftrightarrow \frac{\lambda^{(k_0)}_1}{2}|u_1|^2 + \frac{\lambda^{(k_0)}_2}{2}|u_2|^2  +\mu^1(u_1,u_2) = a(t_{k_0}-t).
\end{align}
If we consider the rescaling
\begin{align}\label{equ-rescaling}
(v_1,v_2) =(\frac{u_1}{\sqrt{t_{k_0}-t}},\frac{u_2}{\sqrt{t_{k_0}-t}}), 
\end{align}
then (\ref{equ-local}) is equivalent to 
\begin{align*}
\frac{\lambda^{(k_0)}_1}{2}|v_1|^2 + \frac{\lambda^{(k_0)}_2}{2}|v_2|^2 = a - \nu_t(v_1,v_2),
\end{align*}
where
\begin{align*}
\nu_t(v_1,v_2) =\frac{\mu^1(\sqrt{t_{k_0}-t}~v_1, \sqrt{t_{k_0}-t}~v_2)}{t_{k_0}-t}.
\end{align*}
By (\ref{equ-limmu}), as $t$ goes to $t_{k_0}$, $\nu_t \colon \mathbb{C}^2 \to \mathbb{R}$ converges $0$ in $C^\infty$-topology on any compact subset on $\mathbb{C}^2$.
Therefore, as $t$ goes to $t_{k_0}$, we have the following convergence 
\begin{align*}
\lim_{t \to t_{k_0}} \frac{1}{\sqrt{t_{k_0}-t}} \big(\mu_{H_{a,b}}^{-1}(c_t) \cap U_{k_0}\big)
= \{ (v_1,v_2) \in \mathbb{C}^2~|~\frac{\lambda^{(k_0)}_1}{2}|v_1|^2 + \frac{\lambda^{(k_0)}_2}{2}|v_2|^2 = a \}.
\end{align*}
Note that $U_{k_0}\cap M^\sigma= \{ (u_1,u_2) \in \mathbb{C}^2~|~u_1,u_2 \in \mathbb{R} \}$ and that the rescaling procedure (\ref{equ-rescaling}) is $H_{a,b}$-equivariant. 
So there exists $\epsilon >0, R>0, C>0$ such that the following $(i)$, $(ii)$ hold:
\\
$(i)$ If $0 < t_{k_0} -t < \epsilon$, then $B(P_{k_0};\sqrt{t_{k_0} -t}~\!R) \cap \mathrm{Im}F_t \ne \emptyset$, where $B(P; r)$ is the geodesic ball of the radius $r>0$ centered at $P$ in $M(\alpha,0)$, 
\\
$(ii)$ $\sup \{|A_t(P)|~|~0 < t_{k_0} -t < \epsilon,~ P \in B(P_{k_0};\sqrt{t_{k_0} -t}~\!R) \cap \mathrm{Im}F_t \}  \le \frac{C}{\sqrt{t_{k_0}-t}}$, where $A_t$ is the second fundamental form of $F_t \colon V_{c_0} \times H_{a,b} \to M(\alpha,0)$. 

Thus we have the following. 
\begin{proposition}\label{prop-sing}
The Lagrangian mean curvature flow $$\{ F_t \colon V_{c_0} \times H_{a,b} \to M(\alpha,0) \}_{t \in [0,t_{k_0})}$$ develops a type I singularity at $P_{k_0}$ when $t$ goes to $t_{k_0}$. The blow-up limit at the singularity is the self-shrinker in Subsection \ref{3-1} with $d=2$, $\lambda_1= a(n-k_0 +1)+b$ and $\lambda_2= -a(n-k_0)-b$. \end{proposition}

As we already mentiond, in \cite{LW1,LW2}, Lee and Wang proved that $\{ L_t \}_{t \in \mathbb{R}}$ in the last paragraph of Subsection  \ref{3-1} forms an eternal solution for Brakke flow. 
Set $L_t^{(a,b)} = \{ ph \in M(\alpha,0)~|~ p \in \mu_{H_{a,b}}^{-1}(c_t) \cap M^\sigma, h \in H_{a,b} \}$ for $t \in \mathbb{R}$. 
By Proposition \ref{prop-sing}, we see that our example $\{ L_t^{(a,b)} \}_{t \in \mathbb{R}}$ becomes closer to $\{ L_t \}_{t \in \mathbb{R}}$, where $d=2$, $\lambda_1= a(n-k_0 +1)+b$ and $\lambda_2= -a(n-k_0)-b$, around $P_{k_0}$ as $t$ goes to $t_{k_0}$. 
So, by modifying the argument of Lee and Wang,  it is not difficult to see that $\{ L_t^{(a,b)} \}_{t \in \mathbb{R}}$ forms an eternal solution for Brakke flow. 
\subsection{The case $a>0,~0 > \frac{b}{a} > -(n+1),~ \frac{b}{a} \not\in \mathbb{Z}$}\label{case2}
Define $t_k \in \mathbb{R}$ by $\iota_{a,b}^*(v_k)=c_0-t_k a w^{a,b}(=c_{t_k})$ for $k=0,1,\dots, n$ as in Subsection \ref{case1}.
Then, in the case $a>0,~0 > \frac{b}{a} > -(n+1),~ \frac{b}{a} \not\in \mathbb{Z}$, there exists $m_0 \in \{ 0, \dots, n \}$ such that 
\begin{align}\label{mzero}
-\infty=t_{-1} < t_0 <t_1 < \dots < t_{m_0} > t_{m_0 +1} > \dots >t_n >t_{n+1}= -\infty.
\end{align}
If $0 > t_{m_0}$, then $\mu_{H_{a,b}}^{-1}(c_0) \cap L = \emptyset$ and we have nothing to discuss. 
So we may assume $t_{m_0} \ge 0$. 
Since we have assumed $\mu_{H_{a,b}}^{-1}(c_0) \cap M(\alpha,0)^G = \emptyset$, there exists $i_0 \in \{0, \dots, m_0 \}$ and $j_0 \in \{ m_0 +1, \dots, n+1 \}$ such that $t_{i_0 -1} < 0 < t_{i_0}$ and $t_{j_0 -1} > 0 > t_{j_0}$ and that $(\iota_{a,b}^*)^{-1}(c_0)$ intersects $\mathrm{int}(l_{i_0})$ and $\mathrm{int}(l_{j_0})$. 
Therefore, in this case, the set $r_{c_0}=(\iota_{a,b}^*)^{-1}(c_0) \cap \mathrm{Im}\mu_G$ is a segment in $\mathfrak{g}^*$. 

By Proposition \ref{prop-fixedpts}, we have the following.
\\
\hspace*{10mm}$(a)$ If  $j_0 -i_0$ is odd,
\\
\hspace*{20mm}$V_{c_0}$ is connected and diffeomorphic to $S^1$.
\\
\hspace*{10mm}$(b)$ If $j_0 -i_0$ is even,
\\
\hspace*{20mm}$V_{c_0}$ consists of two connected components. 
\\
\hspace*{20mm}Each component is diffeomorphic to $S^1$.

If we set $T= \mathrm{min}\{t_{i_0}, t_{j_0 -1} \}$,  then we see that $\mu_{H_{a,b}}^{-1}(c_t) \cap M(\alpha,0)^G = \emptyset$ for each $t \in [0,T)$. 
By Theorem \ref{LMCFALE} we have a Lagrangian mean curvature flow $\{ F_t \colon V_{c_0} \times H_{a,b} \to M(\alpha,0) \}_{t \in [0,T)}$.  
It is impossible to extend the flow continuously to $t \in [0,T+\epsilon)$ for any $\epsilon >0$ as in Subsection \ref{case1}. 
So the Lagrangian mean curvature flow $\{ F_t \colon V_{c_0} \times H_{a,b} \to M(\alpha,0) \}_{t \in [0,T)}$ develops a singularity when $t$ goes to $T$ at $P_{i_0}$ if $T=t_{i_0}$ or at $P_{j_0 -1}$ if $T=t_{j_0 -1}$. 
The structure of the singularity at $P_{k_0}$ is investigated in the same way as in Subsection \ref{case1} and described in Proposition \ref{prop-sing}. 
Here we should note the sign of $\lambda^{(k_0)}_1$ and $\lambda^{(k_0)}_2$ in (\ref{lambda}) as follows. 

\begin{proposition}\label{prop-sign}
Suppose that $a,b \in \mathbb{Z}$ are coprime, $a>0$, $0 > \frac{b}{a} > -(n+1)$ and $\frac{b}{a} \not\in \mathbb{Z}$. 
Then the following holds.
\\
$(1)$ $\lambda^{(m_0)}_1 >0$ and $\lambda^{(m_0)}_2 >0$, where $m_0$ is the same as in (\ref{mzero}). 
\\
$(2)$ If $k_0 \ne m_0$, then $\lambda^{(k_0)}_1 \lambda^{(k_0)}_2 < 0$.
\end{proposition}
{\it Proof.}
For $k_0=0,\dots, n$, we have
\begin{align*}
& \text{$\lambda^{(k_0)}_1 \underset{(\mathrm{resp.}=)}{>}0$ ~~~if and only if~~~ $n+1+\frac{b}{a} \underset{(\mathrm{resp.}=)}{>} k_0$,} \\
& \text{$\lambda^{(k_0)}_2 \underset{(\mathrm{resp.}=)}{>}0$ ~~~if and only if~~~ $k_0 \underset{(\mathrm{resp.}=)}{>} n+\frac{b}{a}$.}
\end{align*}
Since $\frac{b}{a} \not\in \mathbb{Z}$, we have $\lambda^{(k_0)}_1 \ne 0$ and $ \lambda^{(k_0)}_2 \ne 0$. 
So it is enough to show that
\begin{align}\label{slope}
n+1+\frac{b}{a} > m_0 > n+\frac{b}{a}.
\end{align}
In fact, since $\ker \iota_{a,b}^* = \mathrm{span}\{-b p^0 + a p^1 \}$, we see that the slope of the segment $r_{c_0}=(\iota_{a,b}^*)^{-1}(c_0) \cap \mathrm{Im}\mu_G$ is $-\frac{a}{b}$. 
On the other hand, by (\ref{equ-lk}), the slopes of $l_{m_0}$, $l_{m_0 +1}$ are $\frac{1}{n+1-m_0}$, $\frac{1}{n-m_0}$, respectively. 
Then, by  (\ref{mzero}), we have 
\begin{align*}
\frac{1}{n+1-m_0} < -\frac{a}{b} < \frac{1}{n-m_0},
\end{align*}
which is equivalent to (\ref{slope}). 
Thus we finish the proof.
\hfill$\Box$

Set $L_t^{(a,b)} = \{ ph \in M(\alpha,0)~|~ p \in \mu_{H_{a,b}}^{-1}(c_t) \cap M^\sigma, h \in H_{a,b} \}$ for $t \in \mathbb{R}$. 
By Proposition \ref{prop-sing}, it is not difficult to see that $\{ L_t^{(a,b)} \}_{t \in \mathbb{R}}$ forms an eternal solution for Brakke flow as we explained at the end of Subsection \ref{case1}. 
Proposition \ref{prop-sign} implies that $L_t^{(a,b)}$ shrinks to the point $P_{m_0}$ as $t$ goes to $t_{m_0}$ and that $L_t^{(a,b)} = \emptyset$ for $t >t_{m_0}$.
%
%

\begin{picture}(300,315)(0,-15)
\put(0,260){\line(1,0){285}}\put(288,258){\shortstack{$x$}}
\put(40,00){\line(0,1){288}}\put(37,292){\shortstack{$y$}}
\put(43,220){\line(1,0){5}}\put(53,220){\line(1,0){5}}\put(63,220){\line(1,0){5}}\put(73,220){\line(1,0){5}}\put(83,220){\line(1,0){5}}\put(93,220){\line(1,0){5}}
\put(103,220){\line(1,0){5}}\put(113,220){\line(1,0){5}}\put(123,220){\line(1,0){5}}
\put(133,220){\line(1,0){5}}\put(143,220){\line(1,0){5}}\put(153,220){\line(1,0){5}}
\put(163,220){\line(1,0){5}}\put(173,220){\line(1,0){5}}\put(183,220){\line(1,0){5}}
\put(193,220){\line(1,0){5}}\put(203,220){\line(1,0){5}}\put(213,220){\line(1,0){5}}
\put(43,190){\line(1,0){5}}\put(53,190){\line(1,0){5}}\put(63,190){\line(1,0){5}}
\put(73,190){\line(1,0){5}}\put(83,190){\line(1,0){5}}\put(93,190){\line(1,0){5}}
\put(103,190){\line(1,0){5}}\put(113,190){\line(1,0){5}}\put(123,190){\line(1,0){5}}
\put(133,190){\line(1,0){5}}\put(143,190){\line(1,0){5}}\put(153,190){\line(1,0){5}}
\put(43,90){\line(1,0){4}}\put(51,90){\line(1,0){4}}\put(59,90){\line(1,0){4}}
\put(42,60){\line(1,0){2}}\put(46,60){\line(1,0){2}}
\linethickness{1.5pt}
\put(225,220){\line(3,1){60}} \put(253,217){\shortstack{$l_{0}$}}
\put(15,214){\shortstack{$-h_{0}$}}
\put(215,226){\shortstack{$v_{0}$}}
\put(39,220){\line(1,0){2}}
\put(165,190){\line(2,1){60}} \put(193,192){\shortstack{$l_{1}$}}
\put(39,190){\line(1,0){2}}\put(15,184){\shortstack{$-h_{1}$}}
\put(154,195){\shortstack{$v_{1}$}}
\put(39,190){\line(1,0){2}}
\put(149,174){\line(1,1){10}} 
\put(129,154){\line(1,1){10}} 
\put(109,134){\line(1,1){10}} 
\put(89,114){\line(1,1){10}} 
\put(69,94){\line(1,1){10}}  
\put(7,84){\shortstack{$-h_{n-2}$}}
\put(44,95){\shortstack{$v_{n-2}$}}
\put(39,90){\line(1,0){2}}
\put(50,60){\line(1,2){15}} \put(61,67){\shortstack{$l_{n-1}$}}
\put(39,60){\line(1,0){2}} 
\put(7,54){\shortstack{$-h_{n-1}$}}
\put(53,54){\shortstack{$v_{n-1}$}}
\put(39,30){\line(1,0){2}} \put(15,24){\shortstack{$-h_{n}$}}
\put(40,30){\line(1,3){10}} \put(49,38){\shortstack{$l_{n}$}}
\put(43,25){\shortstack{$v_{n}$}}
\put(40,00){\line(0,1){30}} \put(44,9){\shortstack{$l_{n+1}$}}
\put(210,77){\shortstack{$\mathrm{Im} \mu_G$}}
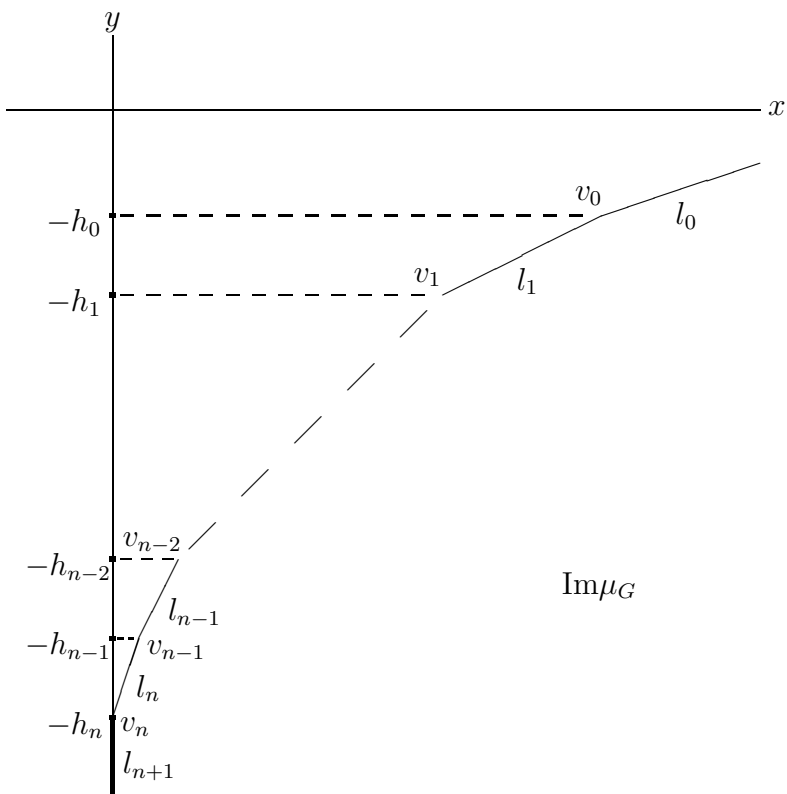
\put(80,-10){\shortstack{Fig. 1.  The image of $\mu_G \colon M(\alpha,0) \to \mathfrak{g}^*$}}
\end{picture}
\end{document}